# RISK BOUNDS FOR STATISTICAL LEARNING

By Pascal Massart and Élodie Nédélec

*Université de Paris-Sud and Université de Paris-Sud*

We propose a general theorem providing upper bounds for the risk of an empirical risk minimizer (ERM).We essentially focus on the binary classification framework. We extend Tsybakov's analysis of the risk of an ERM under margin type conditions by using concentration inequalities for conveniently weighted empirical processes. This allows us to deal with ways of measuring the "size" of a class of classifiers other than entropy with bracketing as in Tsybakov's work. In particular, we derive new risk bounds for the ERM when the classification rules belong to some VC-class under margin conditions and discuss the optimality of these bounds in a minimax sense.

**1. Introduction.** The main results of this paper are obtained within the binary classification framework for which one observes $n$ independent copies $(X_1, Y_1), \ldots, (X_n, Y_n)$ of a pair $(X, Y)$ of random variables, where $X$ takes its values in some measurable space $\mathcal{X}$ and the response variable $Y$ belongs to $\{0, 1\}$. Denoting by $P$ the joint distribution of $(X, Y)$, the quality of a classifier $t$ (i.e., a measurable mapping $t : \mathcal{X} \to \{0, 1\}$) is measured by $P(Y \neq t(X))$. If the distribution $P$ were known, the problem of finding an optimal classifier would be easily solved by considering the Bayes classifier $s^*$ defined for every $x \in \mathcal{X}$ by $s^*(x) = \mathbb{1}_{\eta(x) \geq 1/2}$, where $\eta(x) = P[Y = 1 | X = x]$ denotes the regression function of $Y$ given $X = x$. The Bayes classifier $s^*$ is indeed known to minimize the probability of misclassification $P(Y \neq t(X))$ over the collection of all classifiers. The accuracy of a given classifier $t$ is then measured by its relative loss with respect to the Bayes classifier $\ell(s^*, t) = P(Y \neq t(X)) - P(Y \neq s^*(X))$. The statistical learning problem consists in designing estimators of $s^*$ based on the sample $(X_1, Y_1), \ldots, (X_n, Y_n)$ with as low probability of misclassification as possible. In the sequel we shall use $\ell$ as a loss function and consider the expected risk $\mathbb{E}[\ell(s^*, \hat{s})]$ to analyze the performance of a given estimator $\hat{s}$.

1.1. *Empirical risk minimization.* Given a class of measurable sets $\mathcal{A}$, let $S = \{\mathbb{1}_A, A \in \mathcal{A}\}$ be the corresponding class of classifiers. The empirical risk minimization (ERM) principle consists in taking as an estimator of $s^*$ some minimizer of the empirical criterion

$$t \to \gamma_n(t) = n^{-1} \sum_{i=1}^{n} \mathbb{1}_{Y_i \neq t(X_i)}.$$

Choosing a proper model $S$ among a given list in such a way that simultaneously the bias $\inf_{t \in S} \ell(s^*, t)$ is small enough and the "size" of $S$ is not too large represents the main challenge of model selection procedures. Since the early work of Vapnik and his celebrated book [23], there have been many works on this topic and several attempts to improve on the penalization method of the empirical risk (the structural risk minimization) initially proposed by Vapnik to select among a list of nested models with finite VC-dimensions. Our purpose in this paper is in some sense much less ambitious (although the final goal of our analysis is to design new penalization procedures), and we intend to address the problem of properly identifying what is the benchmark of our estimation problem. More precisely, we just want here to clarify and provide some answers to the following basic questions about the ERM estimators. Assuming first, for the sake of simplicity, that there is no bias, that is, that $s^*$ belongs to $S$:

- What is the order of the expected risk of the ERM on $S$?
- Is it minimax and in what sense?

Of course, since the pioneering work of Vapnik, these questions have been addressed by several authors, but, as we shall see, there are some gaps in the theory. Our aim is to provide some rather complete and general analysis in order to present a unified view allowing a (maybe) better understanding of some already existing results and also to complete the theory by proving some new results.

1.2. *Known risk bounds.* Let us begin with the case where $\mathcal{A}$ is a VC-class.

1.2.1. *Classical bounds for VC-classes.* Recall that if $m_{\mathcal{A}}(N)$ denotes the supremum of $\#\{A \cap C, A \in \mathcal{A}\}$ over the collection of subsets $C$ of $\mathcal{X}$ with cardinality $N$, then $\mathcal{A}$ has the Vapnik–Chervonenkis (VC) property iff $V = \sup\{N : m_{\mathcal{A}}(N) = 2^N\} < \infty$ and $V$ is called the VC-dimension of $\mathcal{A}$. If we denote by $\mathcal{P}(S)$ the set of all joint distributions $P$ such that $s^*$ belongs to $S$ (we must keep in mind that the regression function $\eta$ as well as the Bayes classifier $s^*$ depend on $P$), then (under some convenient measurability condition on $\mathcal{A}$) the following uniform risk bound is available for the ERM $\hat{s}$ (see [13], e.g.):

$$\sup_{P \in \mathcal{P}(S)} \mathbb{E}[\ell(s^*, \hat{s})] \leq \kappa_1 \sqrt{\frac{V}{n}},$$

where $\kappa_1$ denotes some absolute constant. Note that the initial upper bounds on the expected risk for a VC-class found in [23] involved an extra logarithmic factor because they were based on direct combinatorial methods on empirical processes. This factor can be removed (see [13]) by using chaining techniques and the notion of universal entropy, which was introduced in [10] and [19] independently. Furthermore, this upper bound is optimal in the minimax sense since, whenever $2 \leq V \leq n$, one has (see [5])

$$\inf_{\tilde{s}} \sup_{P \in \mathcal{P}(S)} \mathbb{E}[\ell(s^*, \tilde{s})] \geq \kappa_2 \sqrt{\frac{V}{n}},$$

for some absolute positive constant $\kappa_2$, where the infimum is taken over the family of all estimators. Apparently this sounds like the end of the story, but one should realize that this minimax point of view is indeed over-pessimistic. As noted by Vapnik and Chervonenkis themselves in [24], in the (of course, over-optimistic!) situation where $Y = \eta(X)$ almost surely (in this case $P$ is called a zero-error distribution), restricting the set of joint distributions to be with zero-error, the order of magnitude of the minimax lower bound changes drastically, since then one gets $V/n$ instead of $\sqrt{V/n}$. This clearly shows that there is some room for improvment of these global minimax bounds.

### 1.2.2. *Refined bounds for VC-classes.*  Denoting by $\mu$ the marginal distribution of $X$ under $P$, if one takes into account the value of

$$(1) \qquad L(P) = P(Y \neq s^*(X)) = E_\mu[\eta(X) \wedge (1 - \eta(X))],$$

it is possible to get alternative bounds for the risk which can improve on the preceding ones provided that $L(P)$ is small enough [$L(P) = 0$ corresponds to the zero-error case]. The risk bounds found in [5] can be summarized as follows. Given $L_0 \in (0, 1/2)$, if one considers the set $\mathcal{P}_{L_0}(S)$ of distributions $P$ belonging to $\mathcal{P}(S)$ such that $L(P) = L_0$, then (under some measurability condition) for some absolute constant $\kappa_3$, the following upper bound for the risk of the ERM on $S$ is available:

$$(2) \qquad \sup_{P \in \mathcal{P}_{L_0}(S)} \mathbb{E}[\ell(s^*, \hat{s})] \leq \kappa_3 \sqrt{\frac{L_0 V(1 + \log(n/V))}{n}} \qquad \text{if } L_0 \geq \kappa_3(V/n).$$

Moreover, this result is sharp in the minimax sense (up to some logarithmic factor) since, for some absolute positive constant $\kappa_4$, one has

$$(3) \qquad \inf_{\tilde{s}} \sup_{P \in \mathcal{P}_{L_0}(S)} \mathbb{E}[\ell(s^*, \tilde{s})] \geq \kappa_4 \sqrt{\frac{L_0 V}{n}} \qquad \text{if } \kappa_4 L_0 (1 - 2L_0)^2 \geq V/n.$$

We see that (possibly omitting some logarithmic factor) the minimax risk can be of order $V/n$ whenever $L_0$ is of order $V/n$ and that the above bounds offer some kind of interpolation between the zero-error case and the distribution free situations.

However, a careful analysis of the proof of the lower bound (3) shows that the worst distributions are those for which the regression function $\eta$ is allowed to be arbitrarily close to $1/2$. This tends to indicate that maybe some analysis taking into account the way $\eta$ behaves around $1/2$ could be sharper than the preceding one.

1.2.3. *Faster rates under margin conditions.* In [22], Tsybakov attracted attention to rates faster than $1/\sqrt{n}$ that can be achieved by the ERM estimator under a "margin" type condition which is of a different nature from the Devroye and Lugosi condition above, as we shall see below. This condition was first introduced by Mammen and Tsybakov (see [14]) in the related context of discriminant analysis and can be stated as,

$$(4) \qquad \ell(s^*, t) \geq h^\theta \|s^* - t\|_1^\theta \qquad \text{for every } t \in S,$$

where $\| \cdot \|_1$ denotes the $\mathbb{L}_1(\mu)$-norm, $h$ is some positive constant [that we can assume to be smaller than 1 since we can always change $h$ into $h \wedge 1$ without violating (4)] and $\theta \geq 1$. Since $\ell(s^*, t) = E_\mu[|2\eta(X) - 1||s^*(X) - t(X)|]$, we readily see that condition (4) is closely related to the behavior of $\eta(X)$ around $1/2$. In particular, we shall often use in this paper the easily interpretable condition

$$(5) \qquad |2\eta(x) - 1| \geq h \qquad \text{for every } x \in \mathcal{X},$$

which of course implies (4) with $\theta = 1$. Tsybakov uses entropy with bracketing conditions (rather than the VC-condition). In [22], it is shown that, denoting by $H_{[\cdot]}(\varepsilon, S, \mu)$ the $\mathbb{L}_1(\mu)$-entropy with bracketing of $S$ (defined as the logarithm of the minimal number of brackets $[f, g]$ with $\|f - g\|_1 \leq \varepsilon$ which are necessary to cover $S$), if $H_{[\cdot]}(\varepsilon, S, \mu) \ll \varepsilon^{-r}$ for some positive number $r < 1$, then an ERM estimator $\hat{s}$ over $S$ satisfies $\mathbb{E}[\ell(s^*, \hat{s})] = \mathcal{O}(n^{-\theta/(2\theta + r - 1)})$. Hence, Tsybakov's result shows that there is a variety of rates $n^{-\alpha}$ with $1/2 < \alpha < 1$ which can be achieved by an ERM estimator.

### 1.3. *Presentation of our results.*  Our purpose is twofold:

- Providing a general nonasymptotic upper bound for the risk of an ERM which allows to recover Tsybakov's results for classes with integrable entropy with bracketing and also to derive new bounds for VC-classes under margin conditions.
- Focusing on the margin condition (5) and considering $h$ as a free parameter (which may perfectly depend on $n$, for instance), we shall prove minimax lower bounds showing how sharp the preceding upper bounds are.

Even if our upper bounds will cover general margin type conditions [like Tsybakov's condition (4) or even more general than that], we like the idea of focusing on the simpler easy-to-interpret condition (5), which allows comparisons with previous approaches in VC-theory like the one developed in [5]. Let us now state some of the results that we prove in this paper. In order to take into account the margin condition (5) within a minimax approach, we introduce, for every $h \in [0, 1]$, the set $\mathcal{P}(h, S)$ of probability distributions $P$ satisfying the conditions

$$(6) \qquad |2\eta(x) - 1| \geq h \qquad \text{for all } x \in \mathcal{X} \text{ and } s^* \in S$$

(one should keep in mind that $\eta$ as well as $s^*$ depends on $P$, which gives a sense to the definition above). $h = 0$ corresponds to the global minimax approach [one has $\mathcal{P}(0, S) = \mathcal{P}(S)$], while $h = 1$ corresponds to the zero-error case.

### 1.3.1. *The VC-case.*  We assume that $\mathcal{A}$ has finite VC-dimension $V \geq 1$ and we consider an empirical risk minimizer $\hat{s}$ over $S$. Then (at least under some appropriate measurability assumption on the VC-class) we shall prove that, for some absolute positive constant $\kappa$, either

$$\sup_{P \in \mathcal{P}(h,S)} \mathbb{E}[\ell(s^*, \hat{s})] \leq \kappa \sqrt{\frac{V}{n}} \qquad \text{if } h \leq \sqrt{\frac{V}{n}}$$

or

$$(7) \qquad \sup_{P \in \mathcal{P}(h,S)} \mathbb{E}[\ell(s^*, \hat{s})] \leq \kappa \frac{V}{nh}\left(1 + \log\left(\frac{nh^2}{V}\right)\right) \qquad \text{if } h > \sqrt{\frac{V}{n}}.$$

It turns out that, apart from a possible logarithmic factor, this upper bound is optimal in the minimax sense. We indeed show that there exists some absolute positive constant $\kappa'$ such that if $2 \leq V \leq n$,

$$(8) \qquad \inf_{\tilde{s}} \sup_{P \in \mathcal{P}(h,S)} \mathbb{E}[\ell(s^*, \tilde{s})] \geq \kappa'\left[\left(\sqrt{\frac{V}{n}}\right) \wedge \left(\frac{V}{nh}\right)\right].$$

These upper and lower bounds coincide up to the logarithmic factor $1 + \log(nh^2/V)$ and offer some continuous interpolation between the "global" minimax pessimistic bound of order $\sqrt{V/n}$ corresponding to the situation where $h = 0$ (or $h < \sqrt{V/n}$ because then with the margin parameter $h$ being too small, the margin condition has no effect on the order of the minimax risk) and the zero-error case $h = 1$ for which the minimax risk is of order $V/n$ (up to some logarithmic factor).

In order to compare our bounds with those of Devroye and Lugosi [5] recalled above, it is interesting to consider the simple situation where $\eta$ takes on only two values: $(1 - h)/2$ and $(1 + h)/2$. Then by (1), $L(P) = L_0 = (1-h)/2$ so that if, for instance, $h = 1/2$, then $L_0 = 1/4$ and the upper bound given by (7) is of the order of the *square* of the upper bound given by (2). In other words, (2) can be of the same order as in the zero-error case only if $L(P)$ is close enough to zero, while our upper bound is of the same order as in the zero-error case as soon as the margin parameter stays away from 0 (and not only when it is close to 1), which occurs even if $L(P)$ does not tend to zero as $n$ goes to infinity as shown in the preceding elementary example.

We shall also discuss the necessity of the logarithmic factor $1 + \log(nh^2/V)$ in (7). We shall see that the presence of this factor depends on something other than the VC-property. In other words, for some VC-classes, this factor can be removed from the upper bound, while, for some others (which are rich enough in a sense that we shall make precise in Section 3), the minimax lower bound can be refined in order to make this logarithmic factor appear. Quite interestingly, this is, in particular, the case when $\mathcal{A}$ is the class of half-spaces in $\mathbb{R}^d$. We shall indeed prove that in this case, whenever $2 \leq d$, one has, for some positive constant $\kappa''$,

$$
\begin{aligned}
&\inf_{\tilde{s}} \sup_{P \in \mathcal{P}(h,S)} \mathbb{E}[\ell(s^*, \tilde{s})] \\
&\qquad \geq \kappa''(1-h)\left(\frac{d}{nh}\left(1 + \log\left(\frac{nh^2}{d}\right)\right)\right) \qquad \text{if } h \geq \sqrt{d/n}.
\end{aligned}
\tag{9}
$$

We do not know if the factor $1 - h$ in this lower bound can be removed or not but, apart from this factor and up to some absolute positive constant, we can conclude from our study that the minimax risk under the margin condition with parameter $h$ over the class of half-spaces is indeed of order $(d/nh)(1 + \log(nh^2/d))$, provided that $h \geq \sqrt{d/n}$.

1.3.2. *The entropy with bracketing case.* We assume that the entropy with bracketing of $S$ satisfies

$$
H_{[\cdot]}(\varepsilon, S, \mu) \leq K_1 \varepsilon^{-r} \qquad \text{for every } \varepsilon \in (0, 1)
\tag{10}
$$

for some positive number $r < 1$. We can analyze what is the influence on the risk of an ERM $\hat{s}$ of the margin condition (5) by introducing the set $\mathcal{P}(h, S, \mu)$ of distributions $P$ belonging to $\mathcal{P}(h, S)$ with prescribed first marginal distribution $\mu$. Then, for some constant $C_1$ depending only on $K_1$ and $r$, we have

$$\sup_{P \in \mathcal{P}(h, S, \mu)} \mathbb{E}[\ell(s^*, \hat{s})] \leq C_1((nh^{1-r})^{-1/(r+1)} \wedge n^{-1/2}).$$

Moreover, this bound is optimal in the minimax sense, at least if the entropy with bracketing and the $\mathbb{L}_1(\mu)$ metric entropy are of the same order. More precisely, recall that the $\mathbb{L}_1(\mu)$ metric entropy of $S$ denoted by $H_1(\varepsilon, S, \mu)$ is defined as the logarithm of the maximal number of functions $t_1, \ldots, t_N$ belonging to $S$ such that $\|t_i - t_j\|_1 > \varepsilon$ for every $i \neq j$. If (10) holds and if, furthermore, for some positive number $\varepsilon_0 < 1$, one has, for some positive constant $K_2$,

$$(11) \qquad H_1(\varepsilon, S, \mu) \geq K_2 \varepsilon^{-r} \qquad \text{for every } \varepsilon \in (0, \varepsilon_0],$$

then, for some positive constant $C_2$ depending on $K_1, K_2, \varepsilon_0$ and $r$, one has

$$\inf_{\tilde{s}} \sup_{P \in \mathcal{P}(h, S, \mu)} \mathbb{E}[\ell(s^*, \tilde{s})] \geq C_2(1 - h)^{1/(r+1)}((nh^{1-r})^{-1/(r+1)} \wedge n^{-1/2}).$$

In [11, 14] or [6], one can find some explicit examples of classes of subsets of $\mathbb{R}^d$ with smooth boundaries which satisfy both (10) and (11) when $\mu$ is equivalent to the Lebesgue measure on the unit cube.

The paper is organized as follows. In Section 2 we give a general theorem which provides an upper bound for the risk of an ERM via the techniques based on concentration inequalities for weighted empirical processes which were introduced in [15]. The nature of the weight that we are using is absolutely crucial because this is exactly what makes the difference at the end of the day between our upper bounds for VC-classes and those of Devroye and Lugosi [5] which also derive from the analysis of a weighted empirical process but with a different weight. This theorem can be applied to the classification framework, providing the new results described above, but in fact it can also be applied to other frameworks, such as regression with bounded errors. Section 3 is devoted to the minimax lower bounds under margin conditions, while the proofs of all our results are given in Section 4. We have finally postponed to Section 4.2.3 the statements of essentially well known maximal inequalities for empirical processes that we have used all along in the paper.

## 2. A general upper bound for empirical risk minimizers. 
In this section we intend to analyze the behavior of empirical risk minimizers within a framework which is more general than binary classification. Suppose that one observes independent variables $\xi_1, \ldots, \xi_n$ taking their values in some measurable space $\mathcal{Z}$ with common distribution $P$. The two main frameworks that we have in mind are classification and bounded regression. In these cases, for every $i$, the variable $\xi_i = (X_i, Y_i)$ is a copy of a pair of random variables $(X, Y)$, where $X$ takes its values in some measurable space $\mathcal{X}$ and $Y$ is assumed to take its values in $[0, 1]$. In the classification case, the response variable $Y$ is assumed to belong to $\{0, 1\}$. One defines the regression function $\eta$ as $\eta(x) = \mathbb{E}[Y \mid X = x]$ for every $x \in \mathcal{X}$. In the regression case, one is interested in the estimation of $\eta$, while in the classification case one wants to estimate the Bayes classifier $s^*$, defined for every $x \in \mathcal{X}$ by $s^*(x) = \mathbb{1}_{\eta(x) \geq 1/2}$. One of the most commonly used methods to estimate the regression function $\eta$ or the Bayes classifier $s^*$ or, more generally, to estimate a quantity of interest $s$ depending on the unknown distribution $P$, is the so-called empirical risk minimization (according to Vapnik's terminology in [23]). It can be considered a special instance of minimum contrast estimation, which is of course a widely used method in statistics, maximum likelihood estimation being the most celebrated example.

2.1. *Empirical risk minimization.* Basically one considers some set $\mathcal{S}$ which is known to contain $s$. Think of $\mathcal{S}$ as being the set of all measurable functions from $\mathcal{X}$ to $[0, 1]$ in the regression case or to $\{0, 1\}$ in the classification case. Then we consider some *loss (or contrast)* function

$$\gamma \text{ from } \mathcal{S} \times \mathcal{Z} \text{ to } [0, 1], \tag{12}$$

which is well adapted to our problem of estimating $s$ in the sense that the *expected loss* $P[\gamma(t, \cdot)]$ achieves a minimum at the point $s$ when $t$ varies in $\mathcal{S}$. In other words, the *relative expected loss* $\ell$ defined by

$$\ell(s, t) = P[\gamma(t, \cdot) - \gamma(s, \cdot)] \qquad \text{for all } t \in \mathcal{S} \tag{13}$$

is nonnegative. In the regression or the classification case, one can take $\gamma(t, (x, y)) = (y - t(x))^2$ since $\eta$ (resp. $s^*$ ) is indeed the minimizer of $\mathbb{E}[(Y - t(X))^2]$ over the set of measurable functions $t$ taking their values in $[0, 1]$ (resp. $\{0, 1\}$). The heuristics of empirical risk minimization (or minimum contrast estimation) can be described as follows. If one substitutes the empirical loss

$$\gamma_n(t) = P_n[\gamma(t, \cdot)] = \frac{1}{n} \sum_{i=1}^{n} \gamma(t, \xi_i), \tag{14}$$

for its expectation $P[\gamma(t, \cdot)]$ and minimizes $\gamma_n$ on some subset $S$ of $\mathcal{S}$ (that we call a *model*), there is some hope to get a sensible estimator $\hat{s}$ of $s$, at least if $s$ belongs (or is close enough) to model $S$. This estimation method is widely used and has been extensively studied in the asymptotic parametric setting for which one assumes that $S$ is a given parametric model, $s$ belongs to $S$ and $n$ is large.

The purpose of this section is to provide a general nonasymptotic upper bound for the relative expected loss between $\hat{s}$ and $s$.

We introduce the centered empirical process $\overline{\gamma}_n$ defined by

$$(15) \qquad \overline{\gamma}_n(t) = \gamma_n(t) - P[\gamma(t, \cdot)].$$

In addition to the relative expected loss function $\ell$, we shall need another way of measuring the closeness between the elements of $S$ which is directly connected to the variance of the increments of $\overline{\gamma}_n$ and therefore will play an important role in the analysis of the fluctuations of $\overline{\gamma}_n$. Let $d$ be some pseudo-distance on $\mathcal{S} \times \mathcal{S}$ (which may perfectly depend on the unknown distribution $P$) such that

$$(16) \qquad \operatorname{Var}_P[\gamma(t, \cdot) - \gamma(s, \cdot)] \le d^2(s, t) \qquad \text{for every } t \in \mathcal{S}.$$

Of course, we can take $d$ as the pseudo-distance associated with the variance of $\gamma$ itself, but it will be more convenient in applications to take $d$ as a more intrinsic distance. For instance, in the regression or the classification setting it is easy to see that $d$ can be chosen (up to some constant) as the $\mathbb{L}_2(\mu)$ distance, where we recall that $\mu$ denotes the distribution of $X$. Indeed, for classification,

$$|\gamma(t, (x, y)) - \gamma(s^*, (x, y))| = |\mathbb{1}_{y \neq t(x)} - \mathbb{1}_{y \neq s^*(x)}| \le |t(x) - s^*(x)|$$

and, therefore,

$$\operatorname{Var}_P[\gamma(t, \cdot) - \gamma(s^*, \cdot)] \le d^2(s^*, t) \qquad \text{with } d^2(s, t) = E_\mu[(t(X) - s^*(X))^2],$$

while, for regression,

$$[\gamma(t, (x, y)) - \gamma(\eta, (x, y))]^2 = [t(x) - \eta(x)]^2 [2(y - \eta(x)) - t(x) + \eta(x)]^2.$$

Since $E_P[Y - \eta(X) \mid X] = 0$ and $E_P[(Y - \eta(X))^2 \mid X] \le 1/4$, we derive that

$$E_P[[2(Y - \eta(X)) - t(X) + \eta(X)]^2 | X]$$
$$= 4 E_P[(Y - \eta(X))^2 | X] + (-t(X) + \eta(X))^2$$
$$\le 2,$$

and therefore,

$$(17) \qquad E_P[\gamma(t, (X, Y)) - \gamma(\eta, (X, Y))]^2 \le 2 E_\mu (t(X) - \eta(X))^2.$$

Our main result below will crucially depend on two different moduli of uniform continuity: the stochastic modulus of uniform continuity of $\overline{\gamma}_n$ over $S$ with respect to $d$ and the modulus of uniform continuity of $d$ with respect to $\ell$.

The main tool that we shall use is Talagrand's inequality for empirical processes (see [21]) which will allow us to control the oscillations of the empirical process $\overline{\gamma}_n$ by the modulus of uniform continuity of $\overline{\gamma}_n$ in expectation. More precisely, we shall use the following version of it due to Bousquet [4] which has the advantage of providing explicit constants and of dealing with one-sided suprema. If $\mathcal{F}$ is a countable family of measurable functions such that, for some positive constants $v$ and $b$, one has, for every $f \in \mathcal{F}$, $P(f^2) \leq v$ and $\|f\|_\infty \leq b$, then, for every positive $y$, the following inequality holds for $Z = \sup_{f \in \mathcal{F}}(P_n - P)(f)$:

$$(18) \qquad \mathbb{P}\left[ Z - \mathbb{E}[Z] \geq \sqrt{\frac{2(v + 4b\mathbb{E}[Z])y}{n}} + \frac{2by}{3n} \right] \leq e^{-y}.$$

Unlike McDiarmid's inequality (see [18]) which has been widely used in statistical learning theory (see [13]), a concentration inequality like (18) offers the possibility of controlling the empirical process locally. Applying this inequality to some conveniently weighted empirical process will be the key step of the proof of Theorem 2 below.

2.2. *The main theorem.* We need to specify some mild regularity conditions that we shall assume to be verified by the moduli of continuity involved in our result.

DEFINITION 1. We denote by $\mathcal{C}_1$ the class of nondecreasing and continuous functions $\psi$ from $\mathbb{R}_+$ to $\mathbb{R}_+$ such that $x \to \psi(x)/x$ is nonincreasing on $(0, +\infty)$ and $\psi(1) \geq 1$.

Note that if $\psi$ is a nonincreasing continuous and concave function on $\mathbb{R}_+$ with $\psi(0) = 0$ and $\psi(1) \geq 1$, then $\psi$ belongs to $\mathcal{C}_1$. In particular, for the applications that we shall study below, an example of special interest is $\psi(x) = Ax^\alpha$, where $\alpha \in (0, 1]$ and $A \geq 1$.

In order to avoid measurability problems and to use the concentration inequality above, we need to consider some separability condition on $S$. The following one will be convenient.

(M) There exists some countable subset $S'$ of $S$ such that, for every $t \in S$, there exists some sequence $(t_k)$ of elements of $S'$ such that, for every $\xi \in \mathcal{Z}$, $\gamma(t_k, \xi)$ tends to $\gamma(t, \xi)$ as $k$ tends to infinity.

We are now in position to state our upper bound for the relative expected loss of any empirical risk minimizer on some given model $S$. This bound will depend on the bias term $\ell(s, S) = \inf_{t \in S} \ell(s, t)$ and on the fluctuations of the empirical process $\overline{\gamma}_n$ on $S$. As a matter of fact, we shall consider some slightly more general estimators. Namely, given some nonnegative number $\rho$, we consider some $\rho$-empirical risk minimizer, that is, any estimator $\hat{s}$ taking its values in $S$ such that $\gamma_n(\hat{s}) \leq \rho + \inf_{t \in S} \gamma_n(t)$.

THEOREM 2. *Let $\gamma$ be a loss function satisfying* (12) *such that $s$ minimizes $P(\gamma(t, \cdot))$ when $t$ varies in $\mathcal{S}$. Let $\ell$, $\gamma_n$ and $\overline{\gamma}_n$ be defined by* (13), (14) *and* (15) *and consider a pseudo-distance $d$ on $\mathcal{S} \times \mathcal{S}$ satisfying* (16). *Let $\phi$ and $w$ belong to the class of functions $\mathcal{C}_1$ defined above and let $S$ be a subset of $\mathcal{S}$ satisfying the separability condition* (M). *Assume that, on the one hand,*

$$(19) \qquad d(s,t) \leq w(\sqrt{\ell(s,t)}) \qquad \text{for every } t \in \mathcal{S},$$

*and that, on the other hand, one has, for every $u \in S'$,*

$$(20) \qquad \sqrt{n}\mathbb{E}\left[\sup_{t \in S', d(u,t) \leq \sigma}[\overline{\gamma}_n(u) - \overline{\gamma}_n(t)]\right] \leq \phi(\sigma)$$

*for every positive $\sigma$ such that $\phi(\sigma) \leq \sqrt{n}\sigma^2$, where $S'$ is given by assumption* (M). *Let $\varepsilon_*$ be the unique positive solution of the equation*

$$(21) \qquad \sqrt{n}\varepsilon_*^2 = \phi(w(\varepsilon_*)).$$

*Then there exists an absolute constant $\kappa$ such that, for every $y \geq 1$, the following inequality holds:*

$$(22) \qquad \mathbb{P}[\ell(s, \hat{s}) > 2\rho + 2\ell(s, S) + \kappa y \varepsilon_*^2] \leq e^{-y}.$$

*In particular, the following risk bound is available:*

$$\mathbb{E}[\ell(s, \hat{s})] \leq 2(\rho + \ell(s, S) + \kappa \varepsilon_*^2).$$

REMARKS. Let us first give some comments about Theorem 2:

- The absolute constant 2 appearing in (22) has no magic meaning here. It could be replaced by any $C > 1$ at the price of making the constant $\kappa$ depend on $C$.
- One can wonder if an empirical risk minimizer over $S$ exists. Note that condition (M) implies that, for every positive $\rho$, there exists some measurable choice of a $\rho$-empirical risk minimizer since then $\inf_{t \in S'} \gamma_n(t) = \inf_{t \in S} \gamma_n(t)$. If $\rho = 1/n$, for instance, it is clear that, according to (22), such an estimator performs as well as a strict empirical risk minimizer.
- For the computation of $\phi$ satisfying (20), since the supremum appearing in the left-hand side of (20) is extended to the countable set $S'$ and not $S$ itself, it will allow us to restrict ourselves to the case where $S$ is countable.
- It is worth mentioning that, assuming for simplicity that $s \in S$, (22) still holds if we consider the empirical loss $\gamma_n(s) - \gamma_n(\hat{s})$ instead of the expected loss $\ell(s, \hat{s})$. This is indeed a by-product of the proof of Theorem 2 to be found in Section 4.

Even if the main motivation for Theorem 2 is the study of classification, it can also be easily applied to bounded regression. We begin the illustration of Theorem 2 within this framework, which is more elementary than classification since in this case there is a clear connection between the expected loss and the variance of the increments.

2.3. *Application to bounded regression.* In this setting, the regression function $\eta : x \to E_P[Y \mid X = x]$ is the target to be estimated, so that here $s = \eta$. We recall that for this framework we can take $d$ to be the $\mathbb{L}_2(\mu)$ distance times $\sqrt{2}$. The connection between the loss function $\ell$ and $d$ is especially simple in this case. Indeed, $[\gamma(t,(x,y)) - \gamma(\eta,(x,y))] = [-t(x) + \eta(x)][2(y - \eta(x)) - t(x) + \eta(x)]$, so that $E_P[Y - \eta(X) \mid X] = 0$ implies that

$$\ell(\eta, t) = E_P[\gamma(t, (X, Y)) - \gamma(\eta, (X, Y))] = E_\mu(t(X) - \eta(X))^2.$$

Hence, $2\ell(\eta, t) = d^2(\eta, t)$ and in this case the modulus of continuity $w$ can simply be taken as $w(\varepsilon) = \sqrt{2}\varepsilon$. The quadratic risk of an empirical risk minimizer over some model $S$ depends only on the modulus of continuity $\phi$ satisfying (20) and one derives from Theorem 2 that, for some absolute constant $\kappa'$, $\mathbb{E}[d^2(\eta, \hat{s})] \leq 2d^2(\eta, S) + \kappa' \varepsilon_*^2$, where $\varepsilon_*$ is the solution of $\sqrt{n}\varepsilon_*^2 = \phi(\varepsilon_*)$. To be more concrete, let us give an example where this modulus $\phi$ and the bias term $d^2(\eta, S)$ can be evaluated, leading to an upper bound for the minimax risk over some classes of regression functions.

2.3.1. *Binary images.* Following Korostelev and Tsybakov [11], our purpose is to study the particular regression framework for which the variables $X_i$'s are uniformly distributed on $[0,1]^2$ and $\eta(x) = E_P[Y \mid X = x]$ is of the form $\eta(x_1, x_2) = b$ if $x_2 \leq \partial\eta(x_1)$ and $a$ otherwise, where $\partial\eta$ is some measurable map from $[0,1]$ to $[0,1]$ and $0 < a < b < 1$. The function $\partial\eta$ should be understood as the parametrization of a boundary fragment corresponding to some portion $\eta$ of a binary image in the plane ($a$ and $b$ representing the two levels of color which are taken by the image), and restoring this portion of the image from the noisy data $(X_1, Y_1), \ldots, (X_n, Y_n)$ means estimating $\eta$ or, equivalently, $\partial\eta$. Let $\mathcal{G}$ be the set of measurable maps from $[0,1]$ to $[0,1]$. For any $f \in \mathcal{G}$, let us denote by $\chi_f$ the function defined on $[0,1]^2$ by $\chi_f(x_1, x_2) = b$ if $x_2 \leq f(x_1)$ and $a$ otherwise. From this definition, we see that $\chi_{\partial\eta} = \eta$ and, more generally, if we define $\mathcal{S} = \{\chi_f : f \in \mathcal{G}\}$, for every $t \in \mathcal{S}$, we denote by $\partial t$ the element of $\mathcal{G}$ such that $\chi_{\partial t} = t$. It is natural to consider here as an approximate model for $\eta$ a model $S$ of the form $S = \{\chi_f : f \in \partial S\}$, where $\partial S$ denotes some subset of $\mathcal{G}$. Denoting by $\| \cdot \|_1$ (resp. $\| \cdot \|_2$) the Lebesgue $\mathbb{L}_1$-norm (resp. $\mathbb{L}_2$-norm), one has, for every $f, g \in \mathcal{G}$,

$$\|\chi_f - \chi_g\|_1 = (b - a)\|f - g\|_1 \quad \text{and} \quad \|\chi_f - \chi_g\|_2^2 = (b - a)^2\|f - g\|_1$$

or, equivalently, for every $s, t \in \mathcal{S}$,

$$\|s - t\|_1 = (b-a)\|\partial s - \partial t\|_1 \quad \text{and} \quad \|s-t\|_2^2 = (b-a)^2\|\partial s - \partial t\|_1.$$

Given $u = \chi_g \in S$, we have to compute some function $\phi$ satisfying (20) and therefore to majorize $\mathbb{E}[W(\sigma)]$, where $W(\sigma) = \sup_{t \in \mathcal{B}(u,\sigma)} \overline{\gamma}_n(u) - \overline{\gamma}_n(t)$. This can be done using entropy with bracketing arguments. Indeed, let us notice that if $f - \delta \leq f' \leq f + \delta$, then, defining $f_L = \sup(f - \delta, 0)$ and $f_U = \inf(f + \delta, 1)$, the following inequalities hold: $\chi_{f_L} \leq \chi_{f'} \leq \chi_{f_U}$ and $\|\chi_{f_L} - \chi_{f_U}\|_1 \leq 2(b-a)\delta$. This means that, setting $S_\sigma = \{t \in S, d(t, u) \leq \sigma\}$, $\partial S_\rho = \{f \in \partial S, \|f - g\|_1 \leq \rho\}$ and defining $H_\infty(\delta, \rho)$ as the $\mathbb{L}_\infty$ metric entropy for radius $\delta$ of $\partial S_\rho$, one has, for every positive $\varepsilon$,

$$H_{[\cdot]}(\varepsilon, S_\sigma, \mu) \leq H_\infty\left(\frac{\varepsilon}{2(b-a)}, \frac{\sigma^2}{2(b-a)^2}\right).$$

Moreover, if $[t_L, t_U]$ is a bracket with extremities in $\mathcal{S}$ and $\mathbb{L}_1(\mu)$ diameter not larger than $\delta$ and if $t \in [t_L, t_U]$, then

$$y^2 - 2t_U(x)y + t_L^2(x) \leq (y - t(x))^2 \leq y^2 - 2t_L(x)y + t_U^2(x),$$

which implies that $\gamma(\cdot, t)$ belongs to a bracket with $\mathbb{L}_1(P)$-diameter not larger than

$$2E_P\left[(t_U(X) - t_L(X))\left(Y + \frac{t_U(X) + t_L(X)}{2}\right)\right] \leq 2\delta.$$

Hence, if $\mathcal{F} = \{\gamma(\cdot, t), t \in S \text{ and } d(t, u) \leq \sigma\}$, then

$$H_{[\cdot]}(x, \mathcal{F}, P) \leq H_\infty\left(\frac{x}{4(b-a)}, \frac{\sigma^2}{2(b-a)^2}\right)$$

and furthermore, if $d(t, u) \leq \sigma$,

$$\mathbb{E}[|(Y - t(X))^2 - (Y - u(X))^2|] \leq 2\|u - t\|_1 = \frac{2\|u - t\|_2^2}{(b-a)} \leq \frac{\sigma^2}{(b-a)}.$$

Setting

$$\varphi(\sigma) = \int_0^{\sigma/\sqrt{b-a}} \left(H_\infty\left(\frac{x^2}{4(b-a)}, \frac{\sigma^2}{2(b-a)^2}\right)\right)^{1/2} dx,$$

we derive from Lemma A.4 that $\sqrt{n}\mathbb{E}[W(\sigma)] \leq 12\varphi(\sigma)$, provided that

$$(23) \qquad 4\varphi(\sigma) \leq \sqrt{n}\frac{\sigma^2}{(b-a)}.$$

The point now is that, whenever $\partial S$ is part of a linear finite-dimensional subspace of $\mathbb{L}_\infty[0, 1]$, $H_\infty(\delta, \rho)$ is typically bounded by $D[B + \log(\rho/\delta)]$ for some appropriate constants $D$ and $B$. If it is so, then

$$\varphi(\sigma) \le \sqrt{D} \int_0^{\sigma/\sqrt{b-a}} \left( B + \log\left(\frac{2\sigma^2}{x^2(b-a)}\right) \right)^{1/2} dx$$

$$= \frac{\sqrt{2}\sqrt{D}\sigma}{\sqrt{b-a}} \int_0^{1/\sqrt{2}} \sqrt{B + 2|\log(\delta)|}\, d\delta,$$

which implies that, for some absolute constant $\kappa$, $\varphi(\sigma) \le \kappa\sigma\sqrt{(1+B)D/(b-a)}$. The constraint (23) is a fortiori satisfied if $\sigma\sqrt{b-a} \ge 4\kappa\sqrt{(1+B)D/n}$. Hence, if we take $\phi(\sigma) = 12\kappa\sigma\sqrt{(1+B)D/(b-a)}$, assumption (20) is satisfied. To be more concrete, let us consider the example where $\partial S$ is taken to be the set of piecewise constant functions on a regular partition with $D$ pieces on $[0,1]$ with values in $[0,1]$. Then, it is shown in [1] that $H_\infty(\delta, \partial S, \rho) \le D[\log(\rho/\delta)]$ and, therefore, the preceding analysis can be used with $B = 0$. As a matter of fact, this extends to piecewise polynomials with degree not larger than $r$ via some adequate choice of $B$ as a function of $r$, but we just consider the histogram case here to be simple. As a conclusion, Theorem 2 yields in this case for the empirical risk minimizer $\hat{s}$ over $S$

$$\mathbb{E}[\|\partial\eta - \partial\hat{s}\|_1] \le 2\inf_{t \in S} \|\partial\eta - \partial t\|_1 + C\frac{D}{(b-a)^3 n}$$

for some absolute constant $C$. In particular, if $\partial\eta$ satisfies the Hölder smoothness condition $|\partial\eta(x) - \partial\eta(x')| \le L|x - x'|^\alpha$ with $L > 0$ and $\alpha \in (0,1]$, then $\inf_{t \in S} \|\partial\eta - \partial t\|_1 \le LD^{-\alpha}$, leading to

$$\mathbb{E}[\|\partial\eta - \partial\hat{s}\|_1] \le 2LD^{-\alpha} + C\frac{D}{(b-a)^3 n}.$$

Hence, if $\mathcal{H}(L, \alpha)$ denotes the set of functions from $[0,1]$ to $[0,1]$ satisfying the Hölder condition above, an adequate choice of $D$ yields, for some constant $C'$ depending only on $a$ and $b$,

$$\sup_{\partial\eta \in \mathcal{H}(L,\alpha)} \mathbb{E}[\|\partial\eta - \partial\hat{s}\|_1] \le C'\left( L \vee \frac{1}{n} \right)^{1/(\alpha+1)} n^{-\alpha/(1+\alpha)}.$$

As a matter of fact, this upper bound is unimprovable (up to constants) from a minimax point of view (see [11] for the corresponding minimax lower bound).

2.4. *Application to classification.* Our purpose is to apply our main theorem to the classification setting, assuming that the Bayes classifier is the target to be estimated, so that here $s = s^*$. We recall that for this framework we can take $d$ to be the $\mathbb{L}_2(\mu)$-distance (which is also the square root of the $\mathbb{L}_1(\mu)$-distance since we are dealing with $\{0,1\}$-valued functions) and $S = \{\mathbb{1}_A, A \in \mathcal{A}\}$, where $\mathcal{A}$ is some class of measurable sets. Our main task is to compute the moduli of continuity $\phi$ and $w$. In order to evaluate $w$, we need some margin type condition. For instance, we can use Tsybakov's margin condition (4) so that we can also write

$$(24) \qquad \ell(s,t) \geq h^\theta d^{2\theta}(s,t) \qquad \text{for every } t \in S.$$

As quoted in [22], this condition is satisfied if the distribution of $\eta(X)$ is well behaved around $1/2$. If condition (5) holds, then one simply has $\ell(s,t) = E_\mu[|2\eta(X) - 1||s(X) - t(X)|] \geq hd^2(s,t)$, which means that Tsybakov's condition (24) is satisfied with $\theta = 1$. Of course, condition (24) implies that the modulus of continuity $w$ can be taken as

$$(25) \qquad w(\varepsilon) = h^{-1/2}\varepsilon^{1/\theta}.$$

According to the remark following Theorem 2, we shall first assume $S$ to be countable, knowing that our conclusions will remain valid if $S$ is just assumed to satisfy the separability condition (M). In order to evaluate $\phi$, we shall consider two different kinds of assumptions on $S$ which are well known to imply the Donsker property for the class of functions $\{\gamma(t, \cdot), t \in S\}$ and therefore the existence of a modulus $\phi$ which tends to 0 at 0, namely, a Vapnik–Chervonenkis (VC) condition or an entropy with bracketing assumption. Given $u \in S$, in order to bound the expectation of $W(\sigma) = \sup_{d(u,t) \leq \sigma} (-\overline{\gamma}_n(t) + \overline{\gamma}_n(u))$, we shall use the maximal inequalities for empirical processes which are established in the Appendix via slightly different techniques according to the way the "size" of the class $\mathcal{A}$ is measured.

2.4.1. *The VC-case.* We begin with the celebrated VC-condition, which ensures that $\{\gamma(t, \cdot), t \in S\}$ has the Donsker property whatever $P$. So let us assume that $\mathcal{A}$ is a VC-class. One has at least two ways of measuring the "size" of the class $\mathcal{A}$ (or, equivalently, of the class of classifiers $S = \{\mathbb{1}_A, A \in \mathcal{A}\}$):

- The *random combinatorial entropy* defined as

$$H_\mathcal{A} = \log \#\{A \cap \{X_1, \ldots, X_n\}, A \in \mathcal{A}\},$$

which is related to the VC-dimension $V$ of $\mathcal{A}$ via Sauer's lemma (see [13], e.g.) which ensures that

$$H_\mathcal{A} \leq V\left(1 + \log\left(\frac{n}{V}\right)\right)$$

whenever $n \geq V$.

- The Koltchinskii–Pollard notion of *universal metric entropy* defined as follows. For any probability measure $Q$ and every positive $\varepsilon$, let $H_2(\varepsilon, S, Q)$ denote the logarithm of the maximal number of functions $t_1, \ldots, t_N$ belonging to $S$, such that $E_Q(t_i - t_j)^2 > \varepsilon^2$ for every $i \neq j$, and define the universal metric entropy as

$$(26) \qquad H_{\mathrm{univ}}(\varepsilon, S) = \sup_Q H_2(\varepsilon, S, Q),$$

where the supremum is extended to the set of all probability measures on $\mathcal{X}$. The universal metric entropy is related to the VC-dimension via Haussler's bound $H_{\mathrm{univ}}(\varepsilon, \mathcal{A}) \leq \kappa V(1 + \log(\varepsilon^{-1} \vee 1))$, where $\kappa$ denotes some absolute positive constant (see [8]).

The way of expressing $\phi$ by using either the random combinatorial entropy or the universal metric entropy is detailed in the Appendix. Precisely, to use the maximal inequalities proved in the Appendix, we introduce the classes of sets

$$\mathcal{A}_+ = \{\{(x, y) \colon \mathbb{1}_{y \neq t(x)} \leq \mathbb{1}_{y \neq u(x)}\}, t \in S\}$$

and

$$\mathcal{A}_- = \{\{(x, y) \colon \mathbb{1}_{y \neq t(x)} \geq \mathbb{1}_{y \neq u(x)}\}, t \in S\}.$$

Then we define, for every class of sets $\mathcal{B}$ of $\mathcal{X} \times \{0, 1\}$,

$$W_{\mathcal{B}}^+(\sigma) = \sup_{B \in \mathcal{B}, P(B) \leq \sigma^2} (P_n - P)(B) \text{ and } W_{\mathcal{B}}^-(\sigma) = \sup_{B \in \mathcal{B}, P(B) \leq \sigma^2} (P - P_n)(B).$$

Then

$$(27) \qquad \mathbb{E}[W(\sigma)] \leq \mathbb{E}[W_{\mathcal{A}_+}^+(\sigma)] + \mathbb{E}[W_{\mathcal{A}_-}^-(\sigma)]$$

and it remains to control $\mathbb{E}[W_{\mathcal{A}_+}^+(\sigma)]$ and $\mathbb{E}[W_{\mathcal{A}_-}^-(\sigma)]$ via Lemma A.3, which is based either on some direct random combinatorial entropy approach or on some chaining argument and Haussler's bound on the universal entropy recalled above.

More precisely, since the VC-dimensions of $\mathcal{A}_+$ and $\mathcal{A}_-$ are not larger than that of $\mathcal{A}$, and that similarly, the combinatorial entropies of $\mathcal{A}_+$ and $\mathcal{A}_-$ are not larger than the combinatorial entropy of $\mathcal{A}$, denoting by $V$ the VC-dimension of $\mathcal{A}$ (assuming that $V \geq 1$), we derive from (27) and Lemma A.3 that $\sqrt{n}\mathbb{E}[W(\sigma)] \leq \phi(\sigma)$, provided that $\phi(\sigma) \leq \sqrt{n}\sigma^2$, where $\phi$ can be taken either as

$$(28) \qquad \phi(\sigma) = K\sigma\sqrt{(1 \vee \mathbb{E}[H_{\mathcal{A}}])}$$

or as

$$(29) \qquad \phi(\sigma) = K\sigma\sqrt{V(1 + \log(\sigma^{-1} \vee 1))}.$$

In both cases, assumption (20) is satisfied and we can apply Theorem 2 with $w \equiv 1$ or $w$ defined by (25). When $\phi$ is given by (28), the solution $\varepsilon_*$ of equation (21) can be explicitly computed when $w$ is given by (25) or $w \equiv 1$. Hence, the conclusion of Theorem 2 holds with

$$\varepsilon_*^2 = \left(\frac{K^2(1 \vee \mathbb{E}[H_{\mathcal{A}}])}{nh}\right)^{\theta/(2\theta-1)} \wedge \sqrt{\frac{K^2(1 \vee \mathbb{E}[H_{\mathcal{A}}])}{n}}.$$

In the second case, that is, when $\phi$ is given by (29), $w \equiv 1$ implies by (21) that $\varepsilon_*^2 = K\sqrt{V/n}$, while if $w(\varepsilon_*) = h^{-1/2}\varepsilon_*^{1/\theta}$, then

$$\varepsilon_*^2 = K\varepsilon_*^{1/\theta}\sqrt{\frac{V}{nh}}\sqrt{1 + \log((\sqrt{h}\varepsilon_*^{-1/\theta}) \vee 1)}. \tag{30}$$

Since $1 + \log((\sqrt{h}\varepsilon_*^{-1/\theta}) \vee 1) \geq 1$ and $K \geq 1$, we derive from (30) that

$$\varepsilon_*^2 \geq \left(\frac{V}{nh}\right)^{\theta/(2\theta-1)}. \tag{31}$$

Plugging this inequality in the logarithmic factor of (30) yields

$$\varepsilon_*^2 \leq K\varepsilon_*^{1/\theta}\sqrt{\frac{V}{nh}}\sqrt{1 + \frac{1}{2(2\theta-1)}\log\left(\left(\frac{nh^{2\theta}}{V}\right) \vee 1\right)}$$

and, therefore, since $\theta \geq 1$, $\varepsilon_*^2 \leq K\varepsilon_*^{1/\theta}\sqrt{V/(nh)}\sqrt{1 + \log((nh^{2\theta}/V) \vee 1)}$. Hence,

$$\varepsilon_*^2 \leq \left(\frac{K^2V(1 + \log((nh^{2\theta}/V) \vee 1))}{nh}\right)^{\theta/(2\theta-1)}$$
$$\leq K^2\left(\frac{V(1 + \log((nh^{2\theta}/V) \vee 1))}{nh}\right)^{\theta/(2\theta-1)}$$

and, therefore, the conclusion of Theorem 2 holds with

$$\varepsilon_*^2 = K^2\left[\left(\frac{V(1 + \log((nh^{2\theta}/V) \vee 1))}{nh}\right)^{\theta/(2\theta-1)} \wedge \sqrt{\frac{V}{n}}\right].$$

We have a fortiori obtained the following result for the ERM on $S = \{\mathbb{1}_A, A \in \mathcal{A}\}$.

COROLLARY 3. *Assume that $S$ satisfies* (M) *and that $\mathcal{A}$ is a VC-class with dimension $V \geq 1$. There exists an absolute constant $C$ such that if $\hat{s}$ denotes an empirical risk minimizer over $S$ and if $s^*$ belongs to $S$, the following inequality holds:*

$$\mathbb{E}[\ell(s^*, \hat{s})] \leq C\sqrt{\frac{V \wedge (1 \vee \mathbb{E}[H_{\mathcal{A}}])}{n}}. \tag{32}$$

*Moreover, if $\theta \geq 1$ is given and one assumes that the margin condition* (24) *holds with $h \geq (V/n)^{1/2\theta}$, then the following inequalities are also available:*

$$(33) \qquad \mathbb{E}[\ell(s^*, \hat{s})] \leq C \left( \frac{(1 \vee \mathbb{E}[H_{\mathcal{A}}])}{nh} \right)^{\theta/(2\theta-1)}$$

*and*

$$(34) \qquad \mathbb{E}[\ell(s^*, \hat{s})] \leq C \left( \frac{V(1 + \log(nh^{2\theta}/V))}{nh} \right)^{\theta/(2\theta-1)}.$$

Let us comment on these results:

- The risk bound (32) is well known. Our purpose here was just to show how it can be derived from our approach.
- The risk bounds (33) and (34) are new and they perfectly fit with (32) when one considers the borderline case $h = (V/n)^{1/2\theta}$. They look very similar but are not strictly comparable since, roughly speaking, they differ by a logarithmic factor. Indeed, it may happen that $\mathbb{E}[H_{\mathcal{A}}]$ turns out to be of the order of $V$ (without any extra log factor). This is the case when $\mathcal{A}$ is the family of all subsets of a given finite set with cardinality $V$. In such a case, $\mathbb{E}[H_{\mathcal{A}}] \leq V$ and (33) is sharper than (34). On the contrary, for some arbitrary VC-class, if one uses Sauer's bound on $H_{\mathcal{A}}$, that is, $H_{\mathcal{A}} \leq V(1 + \log(n/V))$, the log-factor $1 + \log(n/V)$ is larger than $1 + \log(nh^{2\theta}/V)$ and turns out be too large when $h$ is close to the borderline value $(V/n)^{1/2\theta}$.
- For the sake of simplicity, we have assumed $s^*$ to belong to $S$ in the above statement. Of course, this assumption is not necessary (since our main theorem does not require it). The price to pay if $s^*$ does not belong to $S$ is simply to add $2\ell(s^*, S)$ to the right-hand side of the risk bounds above.

In the next section we shall discuss the optimality of (34) from a minimax point of view in the case where $\theta = 1$, showing that it is essentially unimprovable in that sense.

2.4.2. *Bracketing conditions.* The $\mathbb{L}_1(\mu)$ entropy with bracketing of $S$ is denoted by $H_{[\cdot]}(\delta, S, \mu)$ and has been defined in Section 1. The point is that, setting $\mathcal{F} = \{\gamma(\cdot, t), t \in S$ with $d(u, t) \leq \sigma\}$, one has $H_{[\cdot]}(\delta, \mathcal{F}, P) \leq H_{[\cdot]}(\delta, S, \mu)$. Hence, since we may assume $S$ to be countable (according to the remark after Theorem 2), we derive from (27) and Lemma A.4 in the Appendix that, setting $\varphi(\sigma) = \int_0^\sigma H_{[\cdot]}^{1/2}(x^2, S, \mu) \, dx$, the following inequality is available: $\sqrt{n}\mathbb{E}[W(\sigma)] \leq 12\varphi(\sigma)$, provided that $4\varphi(\sigma) \leq \sigma^2\sqrt{n}$. Hence, we can apply Theorem 2 with $\phi = 12\varphi$, and if we assume Tsybakov's margin condition (24) to be satisfied, then we can also take $w(\varepsilon) = (h^{-1/2}\varepsilon_*^{1/\theta}) \wedge 1$ according to (25) and derive that the conclusions of Theorem 2 hold with $\varepsilon_*$ the solution of the equation $\sqrt{n}\varepsilon_*^2 = \phi((h^{-1/2}\varepsilon_*^{1/\theta}) \wedge 1)$. Moreover, if we assume that condition (10) holds for the entropy with bracketing, then, for some constant $C'$ depending on the constant $K_1$ coming from (10), one has

$$(35) \qquad \varepsilon_*^2 \leq C'[((1-r)^2 n h^{1-r})^{-\theta/(2\theta-1+r)} \wedge (1-r)^{-1} n^{-1/2}].$$

Of course, this conclusion still holds if $S$ is no longer assumed to be countable but fulfills (M). We can alternatively take $T$ to be some $\delta_n$-net [with respect to the $\mathbb{L}_2(\mu)$-distance $d$] of a bigger class $S$ to which the target $s^*$ is assumed to belong. We can still apply Theorem 2 to the empirical risk minimizer over $T$, and since $H_{[\cdot]}(x, T, \mu) \leq H_{[\cdot]}(x, S, \mu)$, we still get the conclusions of Theorem 2 with $\varepsilon_*$ satisfying (35) and $\ell(s^*, T) \leq \delta_n^2$. This means that if $\delta_n$ is conveniently chosen (in a way that $\delta_n$ is of lower order as compared to $\varepsilon_*$), for instance, $\delta_n^2 = n^{-1/(1+r)}$, then, for some constant $C''$ depending only on $K_1$, one has

$$(36) \qquad \mathbb{E}[\ell(s^*, \hat{s})] \leq C''[((1-r)^2 n h^{1-r})^{-\theta/(2\theta-1+r)} \wedge (1-r)^{-1} n^{-1/2}].$$

This means that we have recovered Tsybakov's Theorem 1 in [22] (as a matter of fact, our result is slightly more precise since it also provides the dependence of the risk bound with respect to the margin parameter $h$ and not only on $\theta$ as in Tsybakov's theorem). We refer to [14] for concrete examples of classes of sets with smooth boundaries satisfying (10) when $\mu$ is equivalent to the Lebesgue measure on some compact set of $\mathbb{R}^d$.

## 3. Minimax lower bounds for classification under margin conditions.
We still consider the binary classification framework for which one observes $n$ i.i.d. copies $(X_1, Y_1), \ldots, (X_n, Y_n)$ of a pair of random variables $(X, Y) \in \mathcal{X} \times \{0, 1\}$. The aim is to estimate the Bayes classifier $s^*$. The natural loss function to be considered is

$$\ell(s^*, t) = P(Y \neq t(X)) - P(Y \neq s^*(X)) \geq 0$$

for any classifier $t : \mathcal{X} \to \{0, 1\}$. Our purpose here is to establish lower bounds matching with the upper bounds for the risk of an empirical risk minimizer provided in the preceding section. In particular, we wish to take into account the effect of the margin condition which has been already analyzed for the upper bounds. Toward this aim, we shall use the minimax point of view, but under a convenient margin restriction on the distribution $P$ of the pair $(X, Y)$. Namely, we shall assume that $P$ belongs to the collection of distributions $\mathcal{P}(h, S)$ as defined by (6). If $\mathcal{A}$ denotes the class of sets linked to $S$, that is, $S = \{\mathbb{1}_A, A \in \mathcal{A}\}$, such as for the upper bounds, the way of measuring the size of $\mathcal{A}$ will influence the construction of the lower bounds for the minimax risk

$$(37) \qquad R_n(h, S) = \inf_{\hat{s} \in S} \sup_{P \in \mathcal{P}(h, S)} \mathbb{E}[\ell(s^*, \hat{s})],$$

where the infimum is taken over the set of all estimators based on the $n$-sample $(X_1, Y_1), \ldots, (X_n, Y_n)$ taking their values in $S$. We begin with the case where $\mathcal{A}$ is a VC-class.

3.1. *VC-classes.* We assume $\mathcal{A}$ to be a VC-class of subsets of $\mathcal{X}$ with VC-dimension $V \geq 1$. Some lower bounds for $R_n(h, S)$ are well known in the two extreme cases $h = 0$ and $h = 1$.

If $h = 0$, we are in a (pessimistic) global minimax approach for which there is no margin restriction in fact and the following lower bound can be found in [5]:

$$R_n(0, S) \geq \frac{e^{-8}}{2\sqrt{6}} \sqrt{\frac{V-1}{n}} \tag{38}$$

for every $n \geq 5(V - 1)$.

If $h = 1$, we are in the zero-error case for which $Y = s^*(X)$ and we have at our disposal a lower bound proved by Vapnik and Chervonenkis in [24] (see also [9]),

$$R_n(1, S) \geq \frac{V-1}{4en} \tag{39}$$

for every $n \geq 2 \vee (V - 1)$.

As expected, the order of these lower bounds for the minimax risk is very sensitive to the set of joint distributions over which the supremum is taken. Our purpose is to provide a continuous link between the general case $h = 0$ and the zero-error case $h = 1$. We first prove a lower bound which holds for *any* VC-class $\mathcal{A}$ and then discuss the presence or not of an extra logarithmic factor in the lower bound for some particular examples. As nicely described in [26], there exist several techniques to derive minimax lower bounds in statistics. We shall use two of them below which are based either on Hellinger distance or Kullback–Leibler information computations.

3.1.1. *A general lower bound.* Let $\mathcal{A}$ be some class of measurable subsets of $\mathcal{X}$ and $S$ be the set of classifiers $S = \{\mathbb{1}_A, A \in \mathcal{A}\}$. When $\mathcal{A}$ is an arbitrary VC-class, our lower bound for the minimax risk on $S$ under a margin condition will be obtained via the "Assouad cube" device which involves Hellinger distance computations.

THEOREM 4. *Given $h \in [0, 1]$, we consider the minimax risk $R_n(h, S)$ over the set of distributions $\mathcal{P}(h, S)$ as defined in* (6) *and* (37). *There exists an absolute positive constant $\kappa$ such that, if $\mathcal{A}$ is a VC-class with dimension $V \geq 2$, one has*

$$R_n(h, S) \geq \kappa \left[ \left( \frac{V}{nh} \right) \wedge \sqrt{\frac{V}{n}} \right], \tag{40}$$

*if $n \geq V$.*

The proof of this result will be given in Section 4.2.1. The novelty of this lower bound emerges if $h \geq \sqrt{V/n}$ since then we see that there is indeed an effect of the margin condition as compared to the global bound (38). In particular, for $h = 1$, we recover (39), up to some absolute constant.

Let us now discuss the sharpness of this lower bound by comparing it to the upper bounds derived in the preceding section. Let us consider the very simple example where $\mathcal{A}$ is the collection of all subsets of a given set with cardinality $V$. Then of course $\mathcal{A}$ is a VC-class with dimension $V$, and inequalities (32) and (33) in Corollary 3 ensure that, for some absolute constant $\kappa'$,

$$R_n(h, S) \leq \kappa' \left[ \left( \frac{V}{nh} \right) \wedge \sqrt{\frac{V}{n}} \right].$$

Since (at least if $V \geq 2$) this upper bound coincides with the preceding lower bound up to some absolute constant, this shows that these bounds provide the right order for the minimax risk and therefore cannot be further improved in this case. However, there exist "richer" VC-classes than this one for which the logarithmic factor appearing in (34) is in some sense necessary. This is precisely the purpose of the next section to provide a new combinatorial condition (satisfied by some but not all VC-classes) under which an extra logarithmic factor must appear in the minimax risk.

3.1.2. *A refined lower bound for "rich" VC-classes.* Our purpose is to propose an alternative lower bound for $R_n(h, S)$ when $\mathcal{A}$ is rich enough in a combinatorial sense that we are going to make explicit. Given some integers $D$ and $N$, we introduce the following combinatorial property for the class of sets $\mathcal{A}$:

($\mathrm{A}_{N,D}$) There exist $N$ points $x_1, x_2, \ldots, x_N$ of $\mathcal{X}$ such that the trace of $\mathcal{A}$ on $x = \{x_1, x_2, \ldots, x_N\}$ defined by

$$\mathrm{Tr}(x) = \{A \cap \{x_1, x_2, \ldots, x_N\} : A \in \mathcal{A}\}$$

contains all the subsets of $\{x_1, x_2, \ldots, x_N\}$ with cardinality $D$.

By definition, if $\mathcal{A}$ is a VC-class with dimension $V$, then $\mathcal{A}$ satisfies ($\mathrm{A}_{V,D}$) for all $1 \leq D \leq V$. It is also clear that, given $1 \leq D \leq V$, the VC-class which was analyzed at the end of the preceding section does not satisfy ($\mathrm{A}_{N,D}$) as soon as $N > V$. On the contrary, we shall see below that there are some non-trivial examples of VC-classes which satisfy ($\mathrm{A}_{N,D}$) for arbitrarily large values of $N$ and suitable values of $D$. A convenient information theoretic lemma and combinatorial arguments lead to the following refinement of Theorem 4 that we shall apply to these types of VC-classes. Recall that $S$ denotes the class of classifiers associated with $\mathcal{A}$, that is, $S = \{\mathbb{1}_A; A \in \mathcal{A}\}$.

THEOREM 5. *Given $D \geq 1$, assume that $\mathcal{A}$ satisfies $(A_{N,D})$ for every integer $N$ such that $N \geq 4D$. Given $h \in [0; 1)$, we consider the minimax risk $R_n(h, S)$ over the set of distributions $\mathcal{P}(h, S)$ as defined in (6) and (37). Then, there exists an absolute positive constant $c$ such that*

$$(41) \qquad R_n(h, S) \geq c(1 - h)\frac{D}{nh}\left[1 + \log\left(\frac{nh^2}{D}\right)\right],$$

*provided that*

$$h \geq \sqrt{\frac{D}{n}}.$$

The proof will be given in Section 4.2.2. We intend now to give two explicit examples of VC-classes for which we can apply the preceding lower bound. Given $D \geq 1$, assume $\mathcal{X}$ to be some infinite and countable set and let $\mathcal{A}$ be the collection of all subsets with cardinality $D$ of $\mathcal{X}$. Then $\mathcal{A}$ is a VC-class with dimension $D$ which obviously satisfies property $(A_{N,D})$ for every integer $N \geq D$. Hence,

$$R_n(h, S) \geq c(1 - h)\frac{D}{nh}\left[1 + \log\left(\frac{nh^2}{D}\right)\right],$$

provided that $h \geq \sqrt{D/n}$, and if we compare this lower bound with the upper bound (34), we see that they involve exactly the same logarithmic factor and that they differ by an absolute multiplicative constant times $1 - h$. Thus, apart from this factor $1 - h$, the order of the minimax risk has been identified. As a matter of fact, we do not know how to get rid of this nuisance factor $1 - h$.

This first example could appear to be rather artificial. More interestingly, our result also applies to half-spaces in $\mathbb{R}^d$, for $d \geq 2$. Indeed, a very nice combinatorial geometric result to be found in [7] says that, for every integer $N \geq d + 1$, there exist $N$ distinct points $x_1, x_2, \ldots, x_N$ of $\mathbb{R}^d$ such that the trace of the collection of half-spaces in $\mathbb{R}^d$ on $\{x_1, x_2, \ldots, x_N\}$ contains all the subsets of $\{x_1, x_2, \ldots, x_N\}$ with cardinality $k \leq [d/2]$. This means that the class $\mathcal{A}$ of half-spaces a fortiori satisfies $(A_{N,[d/2]})$, for every integer $N \geq d$. Hence, Theorem 5 applies with $D = [d/2]$ and we get

$$R_n(h, S) \geq \frac{c}{4}(1 - h)\frac{d}{nh}\left[1 + \log\left(\frac{nh^2}{d}\right)\right],$$

provided that $h \geq \sqrt{d/n}$. Furthermore, the VC-dimension of $\mathcal{A}$ is known to be equal to $d + 1$ so that we readily see, as in the preceding example, that the upper bound which derives from (34) coincides with the above lower bound, apart from an absolute constant and possibly the nuisance factor $1 - h$.

The conclusion of the preceding analysis is that the extra logarithmic factor appearing in the upper bound (34) cannot be avoided in general.

3.2. *A lower bound under some purely metric condition.* Our purpose is now to provide a rather general lower bound under some purely metric assumption on $S$ instead of the VC-property. Let $\mathcal{P}(h, S, \mu)$ be the set of distributions $P$ belonging to $\mathcal{P}(h, S)$ with prescribed first marginal distribution $\mu$ and $R_n(h, S, \mu)$ be the corresponding minimax risk

$$R_n(h, S, \mu) = \inf_{\hat{s} \in S} \sup_{P \in \mathcal{P}(h,S,\mu)} \mathbb{E}[\ell(s^*, \hat{s})].$$

The following general result is available.

THEOREM 6. *Let $\mu$ be a probability measure on $\mathcal{X}$ and $S$ be some class of classifiers on $\mathcal{X}$ such that, for some positive constants $K_1, K_2, \varepsilon_0$ and $r$,*

$$K_2 \varepsilon^{-r} \le H_1(\varepsilon, S, \mu) \le K_1 \varepsilon^{-r}$$

*for all $0 < \varepsilon \le \varepsilon_0$, where $H_1(\cdot, S, \mu)$ denotes the $\mathbb{L}_1(\mu)$-metric entropy of $S$. Then, there exists a positive constant $K$ depending on $K_1$, $K_2$, $\varepsilon_0$ and $r$ such that the following lower bound holds:*

$$(42) \quad R_n(h, S, \mu) \ge K(1-h)^{1/(1+r)}[(h^{-(1-r)/(1+r)} n^{-1/(1+r)}) \wedge n^{-1/2}],$$

*whenever $n \ge 2$.*

The proof of this result will be given in Section 4.2.3. If we are in a situation where the $\mathbb{L}_1(\mu)$ metric entropy and the $\mathbb{L}_1(\mu)$ entropy with bracketing are of the same order, we can compare this lower bound with the upper bound (36). More precisely, let us assume that, for some positive constants $K_1, K_2, \varepsilon_0$ and $r < 1$, one has

$$(43) \quad K_2 \varepsilon^{-r} \le H_1(\varepsilon, S, \mu) \le H_{[\cdot]}(\varepsilon, S, \mu) \le K_1 \varepsilon^{-r}$$

for every $\varepsilon \le \varepsilon_0$. Then up to a constant (depending on $K_1, K_2, \varepsilon_0$ and $r$) and the $(1-h)^{1/(r+1)}$ factor, we see that the lower bound (42) and the upper bound (36) coincide. Note that (43) is, in particular, satisfied when $\mathcal{A}$ is a collection of sets with smooth boundaries in various senses as shown in [11, 14] or [6].

## 4. Proofs of the main results.

4.1. *The upper bound*: *proof of Theorem 2.* Since $S$ satisfies (M), we notice that, by dominated convergence, for every $t \in S$, considering the sequence $(t_k)$ provided by condition (M), one has $P(\gamma(\cdot, t_k))$ that tends to $P(\gamma(\cdot, t))$ as $k$ tends to infinity. Hence, $\ell(s, S) = \ell(s, S')$, which implies that there exists some point $\pi(s)$ (which of course may depend on $\varepsilon_*$) such that $\pi(s) \in S'$ and

$$(44) \qquad \ell(s, \pi(s)) \le \ell(s, S) + \varepsilon_*^2.$$

We start from the identity

$$\ell(s,\hat{s}) = \ell(s,\pi(s)) + \gamma_n(\hat{s}) - \gamma_n(\pi(s)) + \overline{\gamma}_n(\pi(s)) - \overline{\gamma}_n(\hat{s}),$$

which, by definition of $\hat{s}$, implies that

$$\ell(s,\hat{s}) \leq \rho + \ell(s,\pi(s)) + \overline{\gamma}_n(\pi(s)) - \overline{\gamma}_n(\hat{s}).$$

Let $x = \sqrt{\kappa' y}\,\varepsilon_*$, where $\kappa'$ is a constant to be chosen later such that $\kappa' \geq 1$ and

$$V_x = \sup_{t \in S} \frac{\overline{\gamma}_n(\pi(s)) - \overline{\gamma}_n(t)}{\ell(s,t) + \varepsilon_*^2 + x^2}.$$

Then,

$$\ell(s,\hat{s}) \leq \rho + \ell(s,\pi(s)) + V_x(\ell(s,\hat{s}) + x^2)$$

and therefore, on the event $V_x < 1/2$, one has

$$\ell(s,\hat{s}) < 2(\rho + \ell(s,\pi(s))) + \varepsilon_*^2 + x^2,$$

yielding

$$(45) \qquad \mathbb{P}[\ell(s,\hat{s}) \geq 2(\rho + \ell(s,S)) + 3\varepsilon_*^2 + x^2] \leq \mathbb{P}[V_x \geq \tfrac{1}{2}].$$

Since $\ell$ is bounded by 1, we may always assume $x$ (and thus $\varepsilon_*$) to be not larger than 1. Assuming that $x \leq 1$, it remains to control the variable $V_x$ via Bousquet's inequality. In order to use Bousquet's inequality, we first remark that, by assumption (M),

$$V_x = \sup_{t \in S'} \frac{\overline{\gamma}_n(\pi(s)) - \overline{\gamma}_n(t)}{\ell(s,t) + \varepsilon_*^2 + x^2},$$

which means that we indeed have to deal with a countably indexed empirical process. Note that the triangle inequality implies via (16), (44) and (19) that

$$(46) \qquad (\mathrm{Var}_P[\gamma(t,\cdot) - \gamma(\pi(s),\cdot)])^{1/2} \leq d(s,t) + d(s,\pi(s))$$

$$\leq 2w(\sqrt{\ell(s,t) + \varepsilon_*^2}).$$

Since $\gamma$ takes its values in $[0,1]$, introducing the function $w_1 = 1 \wedge 2w$, we derive from (46) that

$$\sup_{t \in S} \mathrm{Var}_P\left[\frac{\gamma(t,\cdot) - \gamma(\pi(s),\cdot)}{\ell(s,t) + \varepsilon_*^2 + x^2}\right] \leq \sup_{\varepsilon \geq 0} \frac{w_1^2(\varepsilon)}{(\varepsilon^2 + x^2)^2} \leq \frac{1}{x^2} \sup_{\varepsilon \geq 0} \left(\frac{w_1(\varepsilon)}{\varepsilon \vee x}\right)^2.$$

Now the monotonicity assumptions on $w$ imply that either $w(\varepsilon) \leq w(x)$ if $x \geq \varepsilon$ or $w(\varepsilon)/\varepsilon \leq w(x)/x$ if $x \leq \varepsilon$. Hence, one has in any case $w(\varepsilon)/(\varepsilon \vee x) \leq w(x)/x$, which finally yields

$$\sup_{t \in S} \mathrm{Var}_P\left[\frac{\gamma(t,\cdot) - \gamma(\pi(s),\cdot)}{\ell(s,t) + x^2}\right] \leq \frac{w_1^2(x)}{x^4}.$$

On the other hand, since $\gamma$ takes its values in $[0, 1]$, we have

$$\sup_{t \in S} \left\| \frac{\gamma(t, \cdot) - \gamma(\pi(s), \cdot)}{\ell(s, t) + x^2} \right\|_\infty \leq \frac{1}{x^2}.$$

We can therefore apply (18) with $v = w_1^2(x)x^{-4}$ and $b = x^{-2}$, which gives that, on a set $\Omega_y$ with probability larger than $1 - \exp(-y)$, the inequality

$$(47) \qquad V_x < \mathbb{E}[V_x] + \sqrt{\frac{2(w_1^2(x)x^{-2} + 4\mathbb{E}[V_x])y}{nx^2}} + \frac{y}{nx^2}.$$

Now since $\varepsilon_*$ is assumed to be not larger than 1, one has $w(\varepsilon_*) \geq \varepsilon_*$ and therefore, for every $\sigma \geq w(\varepsilon_*)$, the following inequality derives from the definition of $\varepsilon_*$ by monotonicity:

$$\frac{\phi(\sigma)}{\sigma^2} \leq \frac{\phi(w(\varepsilon_*))}{w^2(\varepsilon_*)} \leq \frac{\phi(w(\varepsilon_*))}{\varepsilon_*^2} = \sqrt{n}.$$

Thus, (20) holds for every $\sigma \geq w(\varepsilon_*)$. In order to control $\mathbb{E}[V_x]$, we intend to use Lemma A.5. For every $t \in S'$, we introduce $a^2(t) = \ell(s, \pi(s)) \vee \ell(s, t)$. Then by (44), $\ell(s, t) \leq a^2(t) \leq \ell(s, t) + \varepsilon_*^2$. Hence, we have, on the one hand, that

$$\mathbb{E}[V_x] \leq \mathbb{E}\left[ \sup_{t \in S'} \frac{\overline{\gamma}_n(\pi(s)) - \overline{\gamma}_n(t)}{a^2(t) + x^2} \right]$$

and, on the other hand, that, for every $\varepsilon \geq \varepsilon_*$,

$$\mathbb{E}\left[ \sup_{t \in S', a(t) \leq \varepsilon} (\overline{\gamma}_n(\pi(s)) - \overline{\gamma}_n(t)) \right] \leq \mathbb{E}\left[ \sup_{t \in S', \ell(s,t) \leq \varepsilon^2} (\overline{\gamma}_n(\pi(s)) - \overline{\gamma}_n(t)) \right].$$

Now by (44) if there exists some $t \in S'$ such that $\ell(s, t) \leq \varepsilon^2$, then $\ell(s, \pi(s)) \leq \varepsilon^2 + \varepsilon_*^2 \leq 2\varepsilon^2$ and therefore, by assumption (19) and monotonicity of $\theta \to w(\theta)/\theta$, $d(\pi(s), t) \leq 2w(\varepsilon\sqrt{2}) \leq 2\sqrt{2}w(\varepsilon)$. Thus, we derive from (20) that, for every $\varepsilon \geq \varepsilon_*$,

$$\mathbb{E}\left[ \sup_{t \in S', \ell(s,t) \leq \varepsilon^2} (\overline{\gamma}_n(\pi(s)) - \overline{\gamma}_n(t)) \right] \leq \phi(2\sqrt{2}w(\varepsilon))$$

and since $\theta \to \phi(2\sqrt{2}w(\theta))/\theta$ is nonincreasing, we can use Lemma A.5 to get

$$\mathbb{E}[V_x] \leq 4\phi(2\sqrt{2}w(x))/(\sqrt{n}x^2),$$

and by monotonicity of $\theta \to \phi(\theta)/\theta$,

$$\mathbb{E}[V_x] \leq 8\sqrt{2}\phi(w(x))/(\sqrt{n}x^2).$$

Thus, using the monotonicity of $\theta \to \phi(w(\theta))/\theta$, and the definition of $\varepsilon_*$, we derive that

$$(48) \qquad \mathbb{E}[V_x] \leq \frac{8\sqrt{2}\phi(w(\varepsilon_*))}{\sqrt{n}x\varepsilon_*} = \frac{8\sqrt{2}\varepsilon_*}{x} \leq \frac{8\sqrt{2}}{\sqrt{\kappa'}y} \leq \frac{8\sqrt{2}}{\sqrt{\kappa'}},$$

provided that $x \geq \varepsilon_*$, which holds since $\kappa' \geq 1$. Now, the monotonicity of $\theta \to w_1(\theta)/\theta$ implies that $x^{-2}w_1^2(x) \leq \varepsilon_*^{-2}w_1^2(\varepsilon_*)$, but since $\phi(\theta)/\theta \geq \phi(1) \geq 1$ for every $\theta \in [0, 1]$, we derive from (21) and the monotonicity of $\phi$ and $\theta \to \phi(\theta)/\theta$ that

$$\frac{w_1^2(\varepsilon_*)}{\varepsilon_*^2} \leq \frac{\phi^2(w_1(\varepsilon_*))}{\varepsilon_*^2} \leq \frac{\phi^2(2w(\varepsilon_*))}{\varepsilon_*^2} \leq 4\frac{\phi^2(w(\varepsilon_*))}{\varepsilon_*^2}$$

and, therefore, $x^{-2}w_1^2(x) \leq 4n\varepsilon_*^2$. Plugging this inequality together with (48) into (47) implies that, on the set $\Omega_y$,

$$V_x < \frac{8\sqrt{2}}{\sqrt{\kappa'}} + \sqrt{\frac{2(4n\varepsilon_*^2 + 32/\sqrt{\kappa'})y}{nx^2}} + \frac{2y}{3nx^2}.$$

It remains to replace $x^2$ by its value $\kappa' y \varepsilon_*^2$ to derive that, on the set $\Omega_y$, the following inequality holds:

$$V_x < \frac{8\sqrt{2}}{\sqrt{\kappa'}} + \sqrt{\frac{8(1 + 4(n\varepsilon_*^2\sqrt{\kappa'})^{-1})}{\kappa'}} + \frac{2}{3\kappa' n\varepsilon_*^2}.$$

Taking into account that $\phi(w(\theta)) \geq \phi(1 \wedge w(\theta)) \geq \theta$ for every $\theta \in [0, 1]$, we deduce from the definition of $\varepsilon_*$ that $n\varepsilon_*^2 \geq 1$ and, therefore, the preceding inequality becomes, on $\Omega_y$,

$$V_x < \frac{8\sqrt{2}}{\sqrt{\kappa'}} + \sqrt{\frac{8(1 + 4/\sqrt{\kappa'})}{\kappa'}} + \frac{2}{3\kappa'}.$$

Hence, choosing $\kappa'$ as a large enough numerical constant warrants that $V_x < 1/2$ on $\Omega_y$ and, therefore, (45) yields

$$\mathbb{P}[\ell(s, \hat{s}) \geq 2(\rho + \ell(s, S)) + x^2 + 3\varepsilon_*^2] \leq \mathbb{P}(\Omega_y^c) \leq e^{-y}.$$

We get the required probability bound (22) by setting $\kappa = \kappa' + 3$. The proof can then be easily completed by integrating the tail bound (22) to derive the required upper bound on the expected risk.

4.2. *Lower bounds.* To prove our various lower bounds, we shall use some particular collections of probability distributions $\{P_t, t \in \mathcal{T}\}$ for the random pair $(X, Y)$ satisfying the margin condition (5). The purpose of the next lemma is to compute the Kullback–Leibler information and the Hellinger distance between pairs of distributions belonging to such a collection.

LEMMA 7. *Let $h \in [0, 1]$, $\mu$ be a probability measure on $\mathcal{X}$ and $\mathcal{T}$ be a collection of classifiers on $\mathcal{X}$. Let $(X, Y)$ be the coordinate mappings on $\mathcal{X} \times \{0, 1\}$, and for every $t \in \mathcal{T}$, define $P_t$ to be the probability distribution on $\mathcal{X} \times \{0, 1\}$ such that, under $P_t$, $X$ has distribution $\mu$ and for every $x \in \mathcal{X}$,*

*Y follows conditionally on $X = x$ a Bernoulli distribution with parameter $\eta_t(x)$. Assume that, for some partition $\mathcal{X} = \mathcal{X}_1 \cup \mathcal{X}_2$, one has $\eta_t(x) = (1 + (2t(x) - 1)h)/2$ for every $x \in \mathcal{X}_1$ and $\eta_t(x) = t(x) = 0$ for every $x \in \mathcal{X}_2$. Denoting by $\|\cdot\|_1$ the $\mathbb{L}_1(\mu)$-norm, for every $s, t \in \mathcal{T}$, the square Hellinger distance between $P_s$ and $P_t$ is given by*

$$\mathcal{H}^2(P_t, P_s) = (1 - \sqrt{1 - h^2})\|t - s\|_1, \tag{49}$$

*while, if $h < 1$, the Kullback–Leibler information between $P_s$ and $P_t$ is given by*

$$\mathcal{K}(P_t, P_s) = h \log\left(\frac{1+h}{1-h}\right)\|t - s\|_1. \tag{50}$$

PROOF.   For every $p \in [0, 1]$, let us denote by $\mathcal{B}(p)$ the Bernoulli distribution with parameter $p$. Then

$$\mathcal{H}^2(\mathcal{B}(p), \mathcal{B}(1-p)) = \mathcal{H}^2(\mathcal{B}(1-p), \mathcal{B}(p)) = 1 - 2\sqrt{p(1-p)},$$

while, if $p \in (0, 1)$,

$$\mathcal{K}(\mathcal{B}(p), \mathcal{B}(1-p)) = \mathcal{K}(\mathcal{B}(1-p), \mathcal{B}(p)) = (1 - 2p)\log\left(\frac{1-p}{p}\right).$$

Setting $p = (1+h)/2$, the point is that, whenever $t(x) \neq s(x)$, either $\eta_t(x) = p$ or $\eta_t(x) = 1 - p$, with $\eta_s(x) = 1 - \eta_t(x)$. Hence,

$$\begin{aligned}
\mathcal{H}^2(P_t, P_s) &= \int_{\mathcal{X}} \mathcal{H}^2(\mathcal{B}(\eta_t(x)), \mathcal{B}(\eta_s(x)))\mathbb{1}_{t(x) \neq s(x)}\, d\mu(x) \\
&= \|t - s\|_1 \mathcal{H}^2(\mathcal{B}(p), \mathcal{B}(1-p)) \\
&= \|t - s\|_1(1 - 2\sqrt{p(1-p)}),
\end{aligned}$$

which leads to (49). Similarly, one has

$$\begin{aligned}
\mathcal{K}(P_t, P_s) &= \int_{\mathcal{X}} \mathcal{K}(\mathcal{B}(\eta_t(x)), \mathcal{B}(\eta_s(x)))\mathbb{1}_{t(x) \neq s(x)}\, d\mu(x) \\
&= \|t - s\|_1 \mathcal{K}(\mathcal{B}(p), \mathcal{B}(1-p)) \\
&= \|t - s\|_1(1 - 2p)\log\left(\frac{1-p}{p}\right),
\end{aligned}$$

which leads to (50), giving the proof of the lemma.   $\square$

Let us now turn to the proof of the lower bound which holds for general VC-classes.

4.2.1. *Proof of Theorem* 4.

Proof.   Since $\mathcal{A}$ is a VC-class with dimension $V$, there exists some set with cardinality $V$ which is shattered by $\mathcal{A}$. Denoting such a set by $\{x_1, x_2, \ldots, x_V\}$, we consider the probability distribution $\mu$ supported by $\{x_1, x_2, \ldots, x_V\}$ and defined by

$$\begin{aligned} \mu(x_i) &= p \qquad \text{for } 1 \le i \le V-1, \\ \mu(x_V) &= 1 - p(V-1), \end{aligned} \tag{51}$$

with $p$ being some nonnegative parameter satisfying $p(V-1) \le 1$ and to be chosen later. Now for every element $b$ of the hyper-cube $\{0,1\}^{V-1}$, let

$$\begin{aligned} \eta_b(x_i) &= \tfrac{1}{2}(1 + (2b_i - 1)h) \qquad \text{for all } 1 \le i \le V-1, \\ \eta_b(x_V) &= 0 \end{aligned}$$

and define $P_b$ as the joint distribution on $\mathcal{X} \times \{0,1\}$ such that, under $P_b$, $X$ has distribution $\mu$ and $Y$ given $X = x_i$ has a Bernoulli distribution with parameter $\eta_b(x_i)$, for every $1 \le i \le V$. The corresponding Bayes classifier $s_b^*$ is given by $s_b^*(x_i) = b_i$ for $1 \le i \le V-1$ and $s_b^*(x_V) = 0$. Since $\{x_1, x_2, \ldots, x_V\}$ is shattered by $\mathcal{A}$, we see that $P_b \in \mathcal{P}(h, S)$ for every $b \in \{0,1\}^{V-1}$. The first step is to relate the minimax risk over $\mathcal{P}(h, S)$ to the minimax risk over the finite subfamily $\{P_b, b \in \{0,1\}^{V-1}\}$ of $\mathcal{P}(h, S)$. Given any classifier $t : \mathcal{X} \to \{0,1\}$, we have $\ell(s_b^*, t) = E_\mu[|2\eta_b(X) - 1| |t(X) - s_b^*(X)|] \ge h\|t - s_b^*\|_1$ and, therefore,

$$R_n(h, S) \ge h \inf_{\hat{s} \in S} \sup_{b \in \{0,1\}^{V-1}} \mathbb{E}_b[\|s_b^* - \hat{s}\|_1].$$

Now given an estimator $\hat{s}$ taking its values in $S$, we can define $\hat{b}$ taking its values in $\{0,1\}^{V-1}$ such that

$$\min_{b' \in \{0,1\}^{V-1}} \|s_{b'}^* - \hat{s}\|_1 = \|s_{\hat{b}}^* - \hat{s}\|_1.$$

Hence, by the triangle inequality,

$$\|s_{\hat{b}}^* - s_b^*\|_1 \le \|s_{\hat{b}}^* - \hat{s}\|_1 + \|s_b^* - \hat{s}\|_1 \le 2\|s_b^* - \hat{s}\|_1,$$

which leads to

$$R_n(h, S) \ge \frac{h}{2} \inf_{\hat{b} \in \{0,1\}^{V-1}} \sup_{b \in \{0,1\}^{V-1}} \mathbb{E}_b[\|s_b^* - s_{\hat{b}}^*\|_1].$$

Moreover, from Lemma 7, for every pair of elements $b, b'$ of the hyper-cube $\{0,1\}^{V-1}$, one has

$$\mathcal{H}^2(P_b, P_{b'}) = (1 - \sqrt{1-h^2})\|s_b^* - s_{b'}^*\|_1 = p(1 - \sqrt{1-h^2})\left(\sum_{i=1}^{V-1} \mathbb{1}_{b_i \ne b_i'}\right).$$

We are now in position to apply Assouad's lemma as stated in [1], for instance. This gives $4R_n(h,S) \geq (V-1)ph[1 - \sqrt{2n\theta}]$, where $\theta = p(1 - \sqrt{1-h^2})$. Since $1 - \sqrt{1-h^2} \leq h^2$, choosing $p = 2/(9nh^2)$ implies that $\sqrt{2n\theta} \leq 2/3$ and, therefore,

$$R_n(h,S) \geq \frac{(V-1)}{54nh},$$

at least if our choice of $p$ satisfies the constraint $p(V-1) \leq 1$, which a fortiori holds whenever $h \geq \sqrt{(V-1)/n}$. It remains to notice that if $h \leq \sqrt{(V-1)/n}$, we can use the preceding construction with $\tilde{h} = \sqrt{(V-1)/n}$ instead of $h$. Of course, the corresponding family $\{P_b, b \in \{0,1\}^{V-1}\}$ is included in $\mathcal{P}(\tilde{h}, S)$, but also in $\mathcal{P}(h, S)$ as well, which means that in this case $R_n(h, S) \geq (V-1)/(54n\tilde{h})$, completing the proof of the result. $\square$

We turn now to the proof of a refined lower bound for classes of sets satisfying $(A_{N,D})$ for every $N \geq 4D$. Fano's lemma is one of the classical tools used to build minimax lower bounds. We would rather use the following very convenient bound for multiple testing due to Birgé (see [2]), which has the advantage of being relevent even when testing only two hypotheses.

LEMMA 8. *Let* $N \geq 1$, $(P_i)_{0 \leq i \leq N}$ *be a family of probability distributions and* $(A_i)_{0 \leq i \leq N}$ *be a family of disjoint events. Let* $a = \min_{0 \leq i \leq N} P_i(A_i)$. *Then, setting* $\overline{\mathcal{K}} = N^{-1} \sum_{i=1}^{N} \mathcal{K}(P_i, P_0)$,

$$\tag{52} a \leq 0.71 \vee \left( \frac{\overline{\mathcal{K}}}{\ln(1+N)} \right).$$

4.2.2. *Proof of Theorem* 5. The basic construction is very similar to the one performed in the proof of Theorem 4 except that, given some integer $N \geq 4D$ to be chosen later, we focus on a particular subset of the hyper-cube $\{0,1\}^N$ instead of the hyper-cube itself. We consider the uniform probability distribution $\mu$ on the set $\{x_1, x_2, \ldots, x_N\}$ provided by assumption $(A_{N,D})$. Moreover, setting

$$\{0,1\}_D^N = \left\{ b \in \{0,1\}^N, \sum_{i=1}^{N} b_i = D \right\},$$

we introduce for every element $b$ of $\{0,1\}_D^N$, $\eta_b(x_i) = \frac{1}{2}(1 + (2b_i - 1)h)$ for all $1 \leq i \leq N$ and define $P_b$ as the joint distribution on $\mathcal{X} \times \{0,1\}$ such that, under $P_b$, $X$ has distribution $\mu$ and $Y$ given $X = x_i$ has a Bernoulli distribution with parameter $\eta_b(x_i)$, for every $1 \leq i \leq N$. The corresponding Bayes classifier $s_b^*$ is given by $s_b^*(x_i) = b_i$ for $1 \leq i \leq N$. Since $\{x_1, x_2, \ldots, x_N\}$ is the set provided by assumption $(A_{N,D})$, we see that $P_b \in \mathcal{P}(h, S)$ for every $b \in \{0,1\}_D^N$. Arguing exactly as in the proof of Theorem 4, we notice that, for any classifier $t$, $\ell(s_b^*, t) \geq h\|t - s_b^*\|_1$, from which we derive that, for any subset $\mathcal{C}$ of $\{0;1\}_D^N$, the following lower bound holds:

$$R_n(h, S) \geq \frac{h}{2} \inf_{\hat{b} \in \mathcal{C}} \sup_{b \in \mathcal{C}} \mathbb{E}_b[\|s_{\hat{b}}^* - s_b^*\|_1].$$

Since for every $b, b' \in \{0,1\}_D^N$,

$$\|s_b^* - s_{b'}^*\|_1 = \frac{1}{N} \left( \sum_{i=1}^N \mathbb{1}_{b_i \neq b'_i} \right) = \frac{1}{N} \delta(b, b'),$$

where $\delta$ denotes Hamming distance on $\{0,1\}_D^N$, one has, for any subset $\mathcal{C}$ of $\{0,1\}_D^N$,

$$(53) \qquad R_n(h, S) \geq \frac{h}{2N} \inf_{\hat{b} \in \mathcal{C}} \sup_{b \in \mathcal{C}} \mathbb{E}_b[\delta(b, \hat{b})],$$

and it remains to construct a set $\mathcal{C}$ with maximal cardinality such that the points of $\mathcal{C}$ are mutually sufficiently distant (w.r.t. the Hamming distance). This can be done thanks to a combinatorial argument due to Birgé and Massart (see [17]). We more precisely use the version of it to be found in [20] and which is more convenient for our needs here. So by Lemma 8 in [20], since $N \geq 4D$, we can choose $\mathcal{C}$ in such a way that

$$(54) \qquad \begin{aligned} &\delta(b, b') > D/2, &&\text{for every } b, b' \text{ in } \mathcal{C} \text{ with } b \neq b', \\ &\log(\#\mathcal{C}) \geq \rho D \log\left(\frac{N}{D}\right), &&\text{where } \rho = 0.233. \end{aligned}$$

For this choice of $\mathcal{C}$, (53) leads to

$$R_n(h, S) \geq \frac{hD}{4N} \inf_{\hat{b} \in \mathcal{C}} \max_{b \in \mathcal{C}} \mathbb{P}_b[b \neq \hat{b}] = \frac{hD}{4N} \inf_{\hat{b} \in \mathcal{C}} \left( 1 - \min_{b \in \mathcal{C}} \mathbb{P}_b[b = \hat{b}] \right).$$

We derive from (52) that, given a point $b_0 \in \mathcal{C}$, for any estimator $\hat{b}$, the following upper bound holds:

$$(55) \qquad \min_{b \in \mathcal{C}} \mathbb{P}_b(\hat{b} = b) \leq \alpha \vee \frac{\overline{\mathcal{K}}}{\log(\#\mathcal{C})},$$

where $\alpha = 0.71$ and

$$\overline{\mathcal{K}} = \frac{1}{\#\mathcal{C} - 1} \sum_{b \in \mathcal{C}, b \neq b_0} \mathcal{K}(P_b^{\otimes n}, P_{b_0}^{\otimes n}) = \frac{n}{\#\mathcal{C} - 1} \sum_{b \in \mathcal{C}, b \neq b_0} \mathcal{K}(P_b, P_{b_0}).$$

For any $b \in \{0,1\}_D^N$, we have $\delta(b, b_0) \leq 2D$, and thanks to Lemma 7, since

$$\|s_b^* - s_{b_0}^*\|_1 = \frac{1}{N} \delta(b, b_0),$$

we derive that

$$\overline{\mathcal{K}} \leq \frac{2Dnh}{N} \log\left(\frac{1+h}{1-h}\right) \leq \frac{4Dh^2n}{(1-h)N}.$$

Combining this inequality with (55) leads to

$$(56) \qquad R_n(h, S) \geq \frac{hD(1-\alpha)}{4N},$$

provided that $4Dh^2n \leq \alpha \log(\#\mathcal{C})(1-h)N$, which via (54) a fortiori holds if

$$(57) \qquad \frac{4nh^2}{(1-h)\alpha\rho} \leq N \log(N/D).$$

It remains to define $N$ in such a way that $4D \leq N$ and (57) holds. Setting

$$N = [a] \qquad \text{with } a = \frac{8nh^2}{(1-h)\rho\alpha(1+\log(nh^2/D))},$$

since $N \leq a$, (56) leads to the desired lower bound (41) with $c = \rho\alpha(1-\alpha)/32$, at least if the constraints $4D \leq N$ and (57) are satisfied. Let us first prove that $a \geq 42D$ (which a fortiori implies that $N \geq 4D$). Indeed, since $x \to x(1+\log(x))^{-1}$ increases on $[1, +\infty)$ and $nh^2 \geq D$, we derive from the definition of $a$ that $a/D \geq 8/(\rho\alpha) \geq 42$ and, therefore, on the one hand, the constraint $N \geq 4D$ is satisfied and, on the other hand,

$$(58) \qquad \frac{N}{a} \geq \frac{a-1}{a} \geq \frac{41}{42}.$$

Now, let us notice that, for $\theta = 21/41$, one has

$$1 + \log(x) \leq x^{1-\theta} \qquad \text{for } x \geq 41$$

[which, by monotonicity, amounts to checking numerically that $1 + \log(x) \leq x^{1-\theta}$ at point $x = 41$] or, equivalently, $\log(x/(1+\log(x))) \geq \theta \log(x)$ for $x \geq 41$. Applying this inequality with

$$x = \frac{4nh^2}{\theta\rho\alpha D} \geq \frac{4}{\theta\rho\alpha} \geq 41$$

leads, by definition of $a$, to

$$\log\left(\frac{a}{2\theta D}\right) \geq \log\left(\frac{x}{1+\log(x)}\right) \geq \theta \log(x) \geq \theta(1+\log(nh^2/D)).$$

Hence, since (58) means that $N \geq a/(2\theta)$, we get

$$N \log(N/D) \geq \frac{a}{2\theta} \log\left(\frac{a}{2\theta D}\right) \geq \frac{a}{2}(1+\log(nh^2/D))$$

and, therefore, (57) holds, completing the proof of Theorem 5.

4.2.3. *Proof of Theorem* 6. For every classifier $t \in S$, we set $\eta_t = (1 + (2t - 1)h)/2$ and define $P_t$ as the joint distribution on $\mathcal{X} \times \{0, 1\}$ such that, under $P_t$, $X$ has distribution $\mu$ and $Y$ given $X = x$ has a Bernoulli distribution with parameter $\eta_t(x)$, for every $x \in \mathcal{X}$. By definition of $P_t$, $t$ is the Bayes classifier related to $P_t$. This shows that the collection $\{P_t, t \in S\}$ is included in $\mathcal{P}(h, S)$. Arguing as in the preceding proofs of lower bounds since, for every classifier $t \in S$, $\ell(s, t) \geq h \|t - s\|_1$, then, for any finite subset $\mathcal{C}$ of $S$, the following lower bound is available:

$$R_n(h, S) \geq \frac{h}{2} \inf_{\hat{s} \in \mathcal{C}} \sup_{s \in \mathcal{C}} \mathbb{E}_s[\|s - \hat{s}\|_1].$$

We use now an argument due to Yang and Barron [25]. Given $\varepsilon > 0$, the idea is to construct an $\varepsilon$-net (i.e., a maximal set of points such that the mutual distances between the elements of this net stay of order $\varepsilon$, less or equal to $2C\varepsilon$, say, for some constant $C > 1$). To do this, we consider an $\varepsilon$-net $\mathcal{C}'$ and a $C\varepsilon$-net $\mathcal{C}''$ of $S$ with respect to the $\mathbb{L}_1(\mu)$-distance. Any point of $\mathcal{C}'$ must belong to some ball with radius $C\varepsilon$ centered at some point of $\mathcal{C}''$. Hence, if $\mathcal{C}$ denotes an intersection of $\mathcal{C}'$ with such a ball with maximal cardinality, one has, for every $t, t' \in \mathcal{C}$ with $t \neq t'$,

$$(59) \qquad \varepsilon \leq \|t - t'\|_1 \leq 2C\varepsilon$$

and

$$(60) \qquad \log(\#\mathcal{C}) \geq H_1(\varepsilon, S, \mu) - H_1(C\varepsilon, S, \mu).$$

Hence, $R_n(h, S) \geq (h\varepsilon/2) \inf_{\hat{s} \in \mathcal{C}} (1 - \inf_{s \in \mathcal{C}} P_s(\hat{s} = s))$ and using again Lemma 8, we derive that $R_n(h, S) \geq (h\varepsilon/2)(1 - \alpha)$, provided that $\overline{\mathcal{K}} \leq \alpha \log(\#\mathcal{C})$, where, given some arbitrary point $t_0$ in $S$,

$$\overline{\mathcal{K}} = \frac{1}{\#\mathcal{C} - 1} \sum_{t \neq t_0} \mathcal{K}(P_t^{\otimes n}, P_{t_0}^{\otimes n})$$

$$= \frac{n}{\#\mathcal{C} - 1} \sum_{t \neq t_0} \mathcal{K}(P_t, P_{t_0}).$$

Thanks to Lemma 7, we know that

$$\overline{\mathcal{K}} \leq 2n \left( \frac{h^2}{1 - h} \right) \sup_{t \in \mathcal{C}} \|t - t_0\|_1$$

$$\leq 8n \left( \frac{h^2}{1 - h} \right) \varepsilon.$$

Now, using our assumption on the behavior of $H_1(\eta, S)$, we easily derive from (60) that properly choosing $C$, for some positive constant $C_1$ (depending on $K_1$, $K_2$ and $r$), one has, for every $\varepsilon \leq \varepsilon_0$, $\log(\#\mathcal{C}) \geq C_1 \varepsilon^{-r}$. Therefore,

$$\frac{\overline{\mathcal{K}}}{\log(\#\mathcal{C})} \leq \frac{8n}{C_1} \left( \frac{h^2}{1 - h} \right) \varepsilon^{1+r}$$

and we can conclude that $R_n(h, S) \geq (h\varepsilon/2)(1 - \alpha)$ whenever

$$\frac{8n}{C_1}\left(\frac{h^2}{1-h}\right)\varepsilon^{1+r} \leq \alpha,$$

that is,

$$\varepsilon \leq \left(\frac{\alpha C_1}{8}\right)^{1/(1+r)} h^{-2/(1+r)}(1-h)^{1/(1+r)}n^{-1/(1+r)}.$$

We may always assume that $\alpha C_1/8 \leq \varepsilon_0^{1+r}$, so that choosing

$$\varepsilon = \left(\frac{\alpha C_1}{8}\right)^{1/(1+r)} h^{-2/(1+r)}(1-h)^{1/(1+r)}n^{-1/(1+r)},$$

the constraint $\varepsilon \leq \varepsilon_0$ is satisfied if we assume that $nh^2 \geq 1$, and we finally get in this case

$$R_n(h, S) \geq \left(\frac{1-\alpha}{2}\right)\left(\frac{\alpha C_1}{8}\right)^{1/(1+r)}(h^{-((1-r)/(1+r))}n^{-(1/(1+r))}(1-h)^{1/(1+r)}).$$

Otherwise, if $nh^2 < 1$, we can always use the preceding lower bound with $\tilde{h} = n^{-1/2}$ instead of $h$, which (at least if $n \geq 2$) leads to $R_n(h, S) \geq C'\sqrt{1/n}$, completing the proof of the lower bound.

## APPENDIX: MAXIMAL INEQUALITIES

Our purpose is here to provide maximal inequalities for set-indexed empirical processes under either the VC-condition or an entropy with bracketing assumption, and also for weighted processes under local conditions. Although these inequalities are essentially well known, we have not always found them explicitly stated in the literature in a way which was satisfactory for our needs. This is the reason why we have decided to remind the reader briefly what these results are and how they can be proved, our feeling being that it could make life easier for a reader who is not familiar with empirical process techniques.

### A.1. Random vectors and Rademacher processes.

A.1.1. *Random vectors.* We recall a simple maximal inequality for random vectors which easily follows from an argument due to Pisier (see [17]). This inequality turns out to be extremely useful for deriving chaining bounds for either sub-Gaussian or empirical processes.

LEMMA A.1. *Let* $(Z_f)_{f \in \mathcal{F}}$ *be a finite family of real-valued random variables. Let* $\psi$ *be a convex and continuously differentiable function on* $[0, b)$ *with* $0 < b \leq +\infty$. *Assume that* $\psi(0) = \psi'(0) = 0$ *and set, for every* $x \geq 0$,

$$\psi^*(x) = \sup_{\lambda \in (0,b)}(\lambda x - \psi(\lambda)).$$

*If for every $\lambda \in (0, b)$ and $f \in \mathcal{F}$, one has*

(A.1) $$\log \mathbb{E}[\exp(\lambda Z_f)] \leq \psi(\lambda),$$

*then, if $N$ denotes the cardinality of $\mathcal{F}$, we have*

$$\mathbb{E}\left[\sup_{f \in \mathcal{F}} Z_f\right] \leq \psi^{*-1}(\log(N)).$$

*In particular, if for some nonnegative number $v$, one has $\psi(\lambda) = \lambda^2 v/2$ for every $\lambda \in (0, +\infty)$, then*

(A.2) $$\mathbb{E}\left(\sup_{f \in \mathcal{F}} Z_f\right) \leq \sqrt{2v \log(N)},$$

*while, if $\psi(\lambda) = \lambda^2 v/(2(1 - c\lambda))$ for every $\lambda \in (0, 1/c)$, one has*

$$\mathbb{E}\left(\sup_{f \in \mathcal{F}} Z_f\right) \leq \sqrt{2v \log(N)} + c \log(N).$$

The two situations where we shall apply this lemma in order to derive chaining bounds are the following:

- $\mathcal{F}$ is a finite subset of $\mathbb{R}^n$ and $Z_f = \sum_{i=1}^{n} \varepsilon_i f_i$, where $(\varepsilon_1, \ldots, \varepsilon_n)$ are independent Rademacher variables. Then, setting $v = \sup_{f \in \mathcal{F}} \sum_{i=1}^{n} f_i^2$, it is well known that (A.1) is satisfied with $\psi(\lambda) = \lambda^2 v/2$ and, therefore,

(A.3) $$\mathbb{E}\left[\sup_{f \in \mathcal{F}} \sum_{i=1}^{n} \varepsilon_i f_i\right] \leq \sqrt{2v \log(N)}.$$

- $\mathcal{F}$ is a finite set of functions $f$ such that $\|f\|_\infty \leq 1$ and $Z_f = \sum_{i=1}^{n} f(\xi_i) - \mathbb{E}[f(\xi_i)]$, where $\xi_1, \ldots, \xi_n$ are independent random variables. Then, setting $v = \sup_{f \in \mathcal{F}} \sum_{i=1}^{n} \mathbb{E}[f^2(\xi_i)]$, as a by-product of the proof of Bernstein's inequality (see [3]), assumption (A.1) is satisfied with $\psi(\lambda) = \lambda^2 v/(2(1 - \lambda/3))$ and, therefore,

(A.4) $$\mathbb{E}\left(\sup_{f \in \mathcal{F}} Z_f\right) \leq \sqrt{2v \log(N)} + \tfrac{1}{3} \log(N).$$

We are now ready to prove a maximal inequality for Rademacher processes which will be useful for analyzing symmetrized empirical processes.

A.1.2. *Rademacher processes.* Let $\mathcal{F}$ be a bounded subset of $\mathbb{R}^n$ equipped with the usual Euclidean norm defined by

$$\|z\|_2^2 = \sum_{i=1}^{n} z_i^2$$

and let, for any positive $\delta$, $H_2(\delta, \mathcal{F})$ denote the logarithm of the maximal number of points $\{f^{(1)}, \ldots, f^{(N)}\}$ belonging to $\mathcal{F}$ such that $\|f^{(j)} - f^{(j')}\|_2^2 > \delta^2$ for every $j \neq j'$. It is easy to derive from the maximal inequality (A.3) the following chaining inequality which is quite standard (see [12]).

LEMMA A.2. *Let $\mathcal{F}$ be a bounded subset of $\mathbb{R}^n$ and $(\varepsilon_1, \ldots, \varepsilon_n)$ be independent Rademacher variables. We consider the Rademacher process $(Z_f)_{f \in \mathcal{F}}$ defined by $Z_f = \sum_{i=1}^n \varepsilon_i f_i$ for every $f \in \mathcal{F}$. Let $\delta$ be such that $\sup_{f \in \mathcal{F}} \|f\|_2 \leq \delta$. Then*

$$(A.5) \qquad \mathbb{E}\left(\sup_{f \in \mathcal{F}} Z_f\right) \leq 3\delta \sum_{j=0}^{\infty} 2^{-j} \sqrt{H_2(2^{-j-1}\delta, \mathcal{F})}.$$

The proof being straightforward, we skip it. The interested reader will find a detailed proof in [16].

We turn now to maximal inequalities for set-indexed empirical processes. The VC-case will be treated via symmetrization by using the preceding bounds for Rademacher processes, while the bracketing case will be studied via a convenient chaining argument.

**A.2. Empirical processes.** Let us first fix some notation. Throughout this section we consider i.i.d. random variables $\xi_1, \ldots, \xi_n$ with values in some measurable space $\mathcal{Z}$ and common distribution $P$. For any $P$-integrable function $f$ on $\mathcal{Z}$, we define $P_n(f) = n^{-1} \sum_{i=1}^n f(\xi_i)$ and $\nu_n(f) = P_n(f) - P(f)$. Given a collection $\mathcal{F}$ of $P$-integrable functions $f$, our purpose is to control the expectation of $\sup_{f \in \mathcal{F}} \nu_n(f)$ or $\sup_{f \in \mathcal{F}} -\nu_n(f)$, when either $\mathcal{F} = \{\mathbb{1}_B, B \in \mathcal{B}\}$ and $\mathcal{B}$ is a VC-class or under an $\mathbb{L}_1$-entropy with bracketing condition on $\mathcal{F}$.

*A.2.1. VC-classes.* In the VC-case the following result is a refinement of what can be found in [15].

LEMMA A.3. *Let $\mathcal{B}$ be a countable VC-class with dimension not larger than $V \geq 1$ and assume that $\sigma > 0$ is such that*

$$P(B) \leq \sigma^2 \qquad \text{for every } B \in \mathcal{B}.$$

*Let*

$$W_{\mathcal{B}}^+ = \sup_{B \in \mathcal{B}} \nu_n(B),$$

$$W_{\mathcal{B}}^- = \sup_{B \in \mathcal{B}} -\nu_n(B)$$

*and*

$$H_{\mathcal{B}} = \log \#\{B \cap \{\xi_1, \ldots, \xi_n\}\}.$$

*Then there exists an absolute constant $K$ such that*

$$(A.6) \qquad \sqrt{n}(\mathbb{E}[W_{\mathcal{B}}^-] \vee \mathbb{E}[W_{\mathcal{B}}^+]) \leq \frac{K}{2}\sigma\sqrt{\mathbb{E}[H_{\mathcal{B}}]},$$

*provided that $\sigma \geq K\sqrt{\mathbb{E}[H_{\mathcal{B}}]/n}$, and*

$$(A.7) \qquad \sqrt{n}(\mathbb{E}[W_{\mathcal{B}}^-] \vee \mathbb{E}[W_{\mathcal{B}}^+]) \leq \frac{K}{2}\sigma\sqrt{V(1 + \log(\sigma^{-1} \vee 1))},$$

*provided that $\sigma \geq K\sqrt{V(1 + |\log \sigma|)/n}$.*

PROOF.   We use the following classical symmetrization device (see, e.g., [12]). Given independent random signs $(\varepsilon_1, \ldots, \varepsilon_n)$, independent of $(\xi_1, \ldots, \xi_n)$, whatever the countable class of functions $\mathcal{F}$, the following inequality holds:

$$(A.8) \qquad \mathbb{E}\left[\sup_{f \in \mathcal{F}}(P_n - P)(f)\right] \leq \frac{2}{n}\mathbb{E}\left[\sup_{f \in \mathcal{F}}\sum_{i=1}^n \varepsilon_i f(\xi_i)\right].$$

Applying this symmetrization inequality to the class $\mathcal{F} = \{\mathbb{1}_B, B \in \mathcal{B}\}$ and the sub-Gaussian inequalities for suprema of Rademacher processes (A.3) or (A.5), setting $\delta_n^2 = [\sup_{B \in \mathcal{B}} P_n(B)] \vee \sigma^2$, we get either

$$(A.9) \qquad \mathbb{E}[W_{\mathcal{B}}^+] \leq 2\sqrt{\frac{2}{n}}\mathbb{E}\sqrt{H_{\mathcal{B}}\delta_n^2}$$

or if $H_{\mathrm{univ}}(\cdot, \mathcal{B})$ denotes the universal entropy of $\mathcal{B}$ as defined in Section 2.4.1,

$$(A.10) \qquad \mathbb{E}[W_{\mathcal{B}}^+] \leq \frac{6}{\sqrt{n}}\mathbb{E}\left[\sqrt{\delta_n^2}\sum_{j=0}^\infty 2^{-j}\sqrt{H_{\mathrm{univ}}(2^{-j-1}\delta_n, \mathcal{B})}\right].$$

Then by the Cauchy–Schwarz inequality, on the one hand, (A.9) becomes

$$(A.11) \qquad \mathbb{E}[W_{\mathcal{B}}^+] \leq 2\sqrt{\frac{2}{n}}\sqrt{\mathbb{E}[H_{\mathcal{B}}]\mathbb{E}[\delta_n^2]},$$

so that

$$\mathbb{E}[W_{\mathcal{B}}^+] \leq 2\sqrt{\frac{2}{n}}\sqrt{\mathbb{E}[H_{\mathcal{B}}](\sigma^2 + \mathbb{E}[W_{\mathcal{B}}^+])},$$

and, on the other hand, since $H_{\mathrm{univ}}(\cdot, \mathcal{B})$ is nonincreasing, we derive from (A.10) that

$$\mathbb{E}[W_{\mathcal{B}}^+] \leq \frac{6}{\sqrt{n}}\sqrt{\mathbb{E}[\delta_n^2]}\sum_{j=0}^\infty 2^{-j}\sqrt{H_{\mathrm{univ}}(2^{-j-1}\sigma, \mathcal{B})},$$

so that, by Haussler's bound (26), one has

$$\mathbb{E}[W_{\mathcal{B}}^+] \leq 6\sqrt{\frac{\kappa V}{n}}\sqrt{\sigma^2 + \mathbb{E}[W_{\mathcal{B}}^+]}\sum_{j=0}^\infty 2^{-j}\sqrt{(j+1)\log(2) + \log(\sigma^{-1} \vee 1) + 1}.$$

Setting either $D = C^2 \mathbb{E}[H_{\mathcal{B}}]$ or $D = C^2 V(1 + \log(\sigma^{-1} \vee 1))$, where $C$ is a conveniently chosen absolute constant $C$ [$C = 2$ in the first case and $C = 6\sqrt{\kappa}(1 + \sqrt{2})$ in the second case], the following inequality holds in both cases:

$$\mathbb{E}[W_{\mathcal{B}}^+] \leq \sqrt{\frac{D}{n}} \sqrt{2[\sigma^2 + \mathbb{E}[W_{\mathcal{B}}^+]]}$$

or, equivalently,

$$\mathbb{E}[W_{\mathcal{B}}^+] \leq \sqrt{\frac{D}{n}} \left[ \sqrt{\frac{D}{n}} + \sqrt{\frac{D}{n} + 2\sigma^2} \right],$$

which, whenever $\sigma \geq 2\sqrt{3}\sqrt{D/n}$, implies that

$$(A.12) \qquad \mathbb{E}[W_{\mathcal{B}}^+] \leq \sqrt{3}\sigma\sqrt{\frac{D}{n}}.$$

The control of $\mathbb{E}[W_{\mathcal{B}}^-]$ is very similar. This time we apply the symmetrization inequality (A.8) to the class $\mathcal{F} = \{-\mathbb{1}_B, B \in \mathcal{B}\}$ and derive by the same arguments as above that

$$\mathbb{E}[W_{\mathcal{B}}^-(\sigma)] \leq \sqrt{\frac{D}{n}} \sqrt{2[\sigma^2 + \mathbb{E}[W_{\mathcal{B}}^+]]}.$$

Hence, provided that $\sigma \geq 2\sqrt{3}\sqrt{D/n}$, (A.12) implies that $\mathbb{E}[W_{\mathcal{B}}^+] \leq \sigma^2/2$ which, in turn, yields $\mathbb{E}[W_{\mathcal{B}}^-] \leq \sqrt{D/n}\sqrt{3\sigma^2}$, completing the proof of the lemma. $\square$

The case of entropy with bracketing can be treated via some direct chaining argument.

A.2.2. *The entropy with bracketing assumption.* We now prove a maximal inequality via a classical chaining argument. Note that the same kind of result would be valid for $\mathbb{L}_2(P)$-entropy with bracketing conditions, but the chaining argument would involve adaptive truncations which are not needed for $\mathbb{L}_1(P)$-entropy with bracketing. Since $\mathbb{L}_1(P)$-entropy with bracketing will suffice for our needs, for the sake of simplicity, we content ourselves with this notion here.

LEMMA A.4. *Let $\mathcal{F}$ be a countable collection of measurable functions such that $0 \leq f \leq 1$ for every $f \in \mathcal{F}$, and let $f_0$ be a measurable function such that $0 \leq f_0 \leq 1$. Let $\delta$ be a positive number such that $P(|f - f_0|) \leq \delta^2$* for every $f \in \mathcal{F}$ and assume $u \to H_{[\cdot]}^{1/2}(u^2, \mathcal{F}, P)$ is integrable at 0. Then, setting

$$\varphi(\delta) = \int_0^\delta H_{[\cdot]}^{1/2}(u^2, \mathcal{F}, P)\,du,$$

the following inequality is available:

$$\sqrt{n}\left(\mathbb{E}\left[\sup_{f \in \mathcal{F}} \nu_n(f_0 - f)\right] \vee \mathbb{E}\left[\sup_{f \in \mathcal{F}} \nu_n(f - f_0)\right]\right) \leq 12\varphi(\delta),$$

provided that $4\varphi(\delta) \leq \delta^2 \sqrt{n}$.

PROOF. We first perform the control of $\mathbb{E}[\sup_{f \in \mathcal{F}} \nu_n(f_0 - f)]$. For any integer $j$, we set $\delta_j = \delta 2^{-j}$ and $H_j = H_{[\cdot]}(\delta_j^2, \mathcal{F}, P)$. By definition of $H_{[\cdot]}(\cdot, \mathcal{F}, P)$, for any integer $j \geq 1$, we can define a mapping $\Pi_j$ from $\mathcal{F}$ to some finite collection of functions such that

$$(A.13) \qquad\qquad \log \#\Pi_j \mathcal{F} \leq H_j$$

and

$$(A.14) \qquad \Pi_j f \leq f \qquad \text{with } P(f - \Pi_j f) \leq \delta_j^2 \text{ for all } f \in \mathcal{F}.$$

For $j = 0$, we choose $\Pi_0$ to be identically equal to $f_0$. For this choice of $\Pi_0$, we still have

$$(A.15) \qquad P(|f - \Pi_0 f|) = P(|f - f_0|) \leq \delta_0^2 = \delta$$

for every $f \in \mathcal{F}$. Furthermore, since we may always assume that the extremities of the brackets used to cover $\mathcal{F}$ take their values in $[0, 1]$, we also have for every integer $j$ that

$$0 \leq \Pi_j f \leq 1.$$

Noticing that since $u \to H_{[\cdot]}(u^2, \mathcal{F}, P)$ is nonincreasing,

$$H_1 \leq \delta_1^{-2} \varphi^2(\delta),$$

and under the condition $4\varphi(\delta) \leq \delta^2 \sqrt{n}$, one has $H_1 \leq \delta_1^2 n$. Thus, since $j \to H_j \delta_j^{-2}$ increases to infinity, the set $\{j \geq 0 : H_j \leq \delta_j^2 n\}$ is a nonvoid interval of the form

$$\{j \geq 0 : H_j \leq \delta_j^2 n\} = [0, J],$$

with $J \geq 1$. For every $f \in \mathcal{F}$, starting from the decomposition

$$-\nu_n(f) = \sum_{j=0}^{J-1} \nu_n(\Pi_j f) - \nu_n(\Pi_{j+1} f) + \nu_n(\Pi_J f) - \nu_n(f),$$

we derive, since $\Pi_J(f) \leq f$ and $P(f - \Pi_J(f)) \leq \delta_J^2$, that

$$-\nu_n(f) \leq \sum_{j=0}^{J-1} \nu_n(\Pi_j f) - \nu_n(\Pi_{j+1} f) + \delta_J^2$$

and, therefore,

(A.16)
$$\mathbb{E}\left[\sup_{f \in \mathcal{F}}[-\nu_n(f)]\right]$$
$$\leq \sum_{j=0}^{J-1} \mathbb{E}\left[\sup_{f \in \mathcal{F}}[\nu_n(\Pi_j f) - \nu_n(\Pi_{j+1} f)]\right] + \delta_J^2.$$

Now, it follows from (A.14) and (A.15) that, for every integer $j$ and every $f \in \mathcal{F}$, one has

$$P[|\Pi_j f - \Pi_{j+1} f|] \leq \delta_j^2 + \delta_{j+1}^2 = 5\delta_{j+1}^2$$

and, therefore, since $|\Pi_j f - \Pi_{j+1} f| \leq 1$,

$$P[|\Pi_j f - \Pi_{j+1} f|^2] \leq 5\delta_{j+1}^2.$$

Moreover, (A.13) ensures that the number of functions of the form $\Pi_j f - \Pi_{j+1} f$ when $f$ varies in $\mathcal{F}$ is not larger than $\exp(H_j + H_{j+1}) \leq \exp(2H_{j+1})$. Hence, we derive from (A.4) that

$$\sqrt{n}\mathbb{E}\left[\sup_{f \in \mathcal{F}}[\nu_n(\Pi_j f) - \nu_n(\Pi_{j+1} f)]\right]$$
$$\leq 2\left[\delta_{j+1}\sqrt{5H_{j+1}} + \frac{1}{3\sqrt{n}}H_{j+1}\right]$$

and (A.16) becomes

(A.17)
$$\sqrt{n}\mathbb{E}\left[\sup_{f \in \mathcal{F}}[-\nu_n(f)]\right]$$
$$\leq 2\sum_{j=1}^{J}\left[\delta_j\sqrt{5H_j} + \frac{1}{3\sqrt{n}}H_j\right] + 4\sqrt{n}\delta_{J+1}^2.$$

It follows from the definition of $J$ that, on the one hand, for every $j \leq J$,

$$\frac{1}{3\sqrt{n}}H_j \leq \frac{1}{3}\delta_j\sqrt{H_j}$$

and, on the other hand,

$$4\sqrt{n}\delta_{J+1}^2 \leq 4\delta_{J+1}\sqrt{H_{J+1}}.$$

Hence, plugging these inequalities in (A.17) yields

$$\sqrt{n}\mathbb{E}\Big[\sup_{f\in\mathcal{F}}[-\nu_n(f)]\Big] \leq 6\sum_{j=1}^{J+1}[\delta_j\sqrt{H_j}],$$

and the result follows. The control of $\mathbb{E}[\sup_{f\in\mathcal{F}}\nu_n(f-f_0)]$ can be performed analogously, changing lower into upper approximations in the dyadic approximation scheme described above. $\square$

**A.3. A maximal inequality for weighted processes.** The following inequality is more or less classical and well known. We present a (short) proof for the sake of completeness. Note that in the statement and the proof of Lemma A.5 below we use the convention that $\sup_{t\in A} g(t) = 0$ whenever $A$ is the empty set.

LEMMA A.5. *Let $S$ be a countable set, $u \in S$ and $a : S \to \mathbb{R}_+$ such that $a(u) = \inf_{t\in S} a(t)$. Let $Z$ be a process indexed by $S$ and assume that the nonnegative random variable $\sup_{t\in\mathcal{B}(\varepsilon)}[Z(u)-Z(t)]$ has finite expectation for any positive number $\varepsilon$, where $\mathcal{B}(\varepsilon) = \{t \in S, a(t) \leq \varepsilon\}$. Let $\psi$ be a nonnegative function on $\mathbb{R}_+$ such that $\psi(x)/x$ is nonincreasing on $\mathbb{R}_+$ and satisfies for some positive number $\varepsilon_*$*

$$\mathbb{E}\Big[\sup_{t\in\mathcal{B}(\varepsilon)}[Z(u)-Z(t)]\Big] \leq \psi(\varepsilon) \qquad \text{for any } \varepsilon \geq \varepsilon_*.$$

*Then, one has, for any positive number $x \geq \varepsilon_*$,*

$$\mathbb{E}\Big[\sup_{t\in S}\Big[\frac{Z(u)-Z(t)}{a^2(t)+x^2}\Big]\Big] \leq 4x^{-2}\psi(x).$$

PROOF.    Let us introduce for any integer $j$

$$\mathcal{C}_j = \{t \in S, r^j x < a(t) \leq r^{j+1}x\},$$

with $r > 1$ to be chosen later. Then $\{\mathcal{B}_u(x), \{\mathcal{C}_j\}_{j\geq 0}\}$ is a partition of $S$ and, therefore,

$$\sup_{t\in S}\Big[\frac{Z(u)-Z(t)}{a^2(t)+x^2}\Big] \leq \sup_{t\in\mathcal{B}_u(x)}\Big[\frac{(Z(u)-Z(t))_+}{a^2(t)+x^2}\Big]$$
$$+ \sum_{j\geq 0}\sup_{t\in\mathcal{C}_j}\Big[\frac{(Z(u)-Z(t))_+}{a^2(t)+x^2}\Big],$$

which, in turn, implies that

$$x^2\sup_{t\in S}\Big[\frac{Z(u)-Z(t)}{a^2(t)+x^2}\Big]$$

$$(A.18) \qquad \leq \sup_{t \in \mathcal{B}_u(x)} (Z(u) - Z(t))_+$$

$$+ \sum_{j \geq 0} (1 + r^{2j})^{-1} \sup_{t \in \mathcal{B}_u(r^{j+1}x)} (Z(u) - Z(t))_+.$$

Since $a(u) = \inf_{t \in S} a(t)$, one has $u \in \mathcal{B}_u(r^k x)$ for every integer $k$ for which $\mathcal{B}_u(r^k x)$ is nonempty and, therefore,

$$\sup_{t \in \mathcal{B}_u(r^k x)} (Z(u) - Z(t))_+ = \sup_{t \in \mathcal{B}_u(r^k x)} (Z(u) - Z(t)).$$

Hence, taking the expectation in (A.18) yields

$$x^2 \mathbb{E}\left[\sup_{t \in S}\left[\frac{Z(u) - Z(t)}{a^2(t) + x^2}\right]\right] \leq \psi(x) + \sum_{j \geq 0} (1 + r^{2j})^{-1} \psi(r^{j+1}x).$$

Now by our monotonicity assumption, $\psi(r^{j+1}x) \leq r^{j+1}\psi(x)$, and thus

$$x^2 \mathbb{E}\left[\sup_{t \in S}\left[\frac{Z(u) - Z(t)}{a^2(t) + x^2}\right]\right] \leq \psi(x)\left[1 + r\sum_{j \geq 0} r^j (1 + r^{2j})^{-1}\right]$$

$$\leq \psi(x)\left[1 + r\left(\frac{1}{2} + \sum_{j \geq 1} r^{-j}\right)\right]$$

$$\leq \psi(x)\left[1 + r\left(\frac{1}{2} + \frac{1}{r-1}\right)\right],$$

and the result follows by choosing $r = 1 + \sqrt{2}$.   $\square$

**Acknowledgments.** We would like to thank Stéphane Boucheron for having encouraged us to write this paper and also for numerous and fruitful discussions on the topic of statistical learning. We also thank the referees for their helpful suggestions.

## REFERENCES

[1] BARRON, A. R., BIRGÉ, L. and MASSART, P. (1999). Risk bounds for model selection via penalization. *Probab. Theory Related Fields* **113** 301–413. MR1679028

[2] BIRGÉ, L. (2005). A new lower bound for multiple hypothesis testing. *IEEE Trans. Inform. Theory* **51** 1611–1615. MR2241522

[3] BIRGÉ, L. and MASSART, P. (1998). Minimum contrast estimators on sieves: Exponential bounds and rates of convergence. *Bernoulli* **4** 329–375. MR1653272

[4] BOUSQUET, O. (2002). A Bennett concentration inequality and its application to suprema of empirical processes. *C. R. Math. Acad. Sci. Paris* **334** 495–500. MR1890640

[5] DEVROYE, L. and LUGOSI, G. (1995). Lower bounds in pattern recognition and learning. *Pattern Recognition* **28** 1011–1018.

[6] DUDLEY, R. M. (1999). *Uniform Central Limit Theorems.* Cambridge Univ. Press. MR1720712

[7] EDELSBRUNNER, H. (1987). *Algorithms in Combinatorial Geometry.* Springer, Berlin. MR0904271

[8] HAUSSLER, D. (1995). Sphere packing numbers for subsets of the Boolean *n*-cube with bounded Vapnik–Chervonenkis dimension. *J. Combin. Theory Ser. A* **69** 217–232. MR1313896

[9] HAUSSLER, D., LITTLESTONE, N. and WARMUTH, M. (1994). Predicting $\{0, 1\}$-functions on randomly drawn points. *Inform. and Comput.* **115** 248–292. MR1304811

[10] KOLTCHINSKII, V. I. (1981). On the central limit theorem for empirical measures. *Theor. Probab. Math. Statist.* **24** 71–82. MR0628431

[11] KOROSTELEV, A. P. and TSYBAKOV, A. B. (1993). *Minimax Theory of Image Reconstruction. Lecture Notes in Statist.* **82**. Springer, New York. MR1226450

[12] LEDOUX, M. and TALAGRAND, M. (1991). *Probability in Banach Spaces. Isoperimetry and Processes.* Springer, Berlin. MR1102015

[13] LUGOSI, G. (2002). Pattern classification and learning theory. In *Principles of Nonparametric Learning* (L. Györfi, ed.) 1–56. Springer, Vienna. MR1987656

[14] MAMMEN, E. and TSYBAKOV, A. B. (1999). Smooth discrimination analysis. *Ann. Statist.* **27** 1808–1829. MR1765618

[15] MASSART, P. (2000). Some applications of concentration inequalities to statistics. Probability theory. *Ann. Fac. Sci. Toulouse Math. (6)* **9** 245–303. MR1813803

[16] MASSART, P. (2006). Concentration inequalities and model selection. *Lectures on Probability Theory and Statistics. Ecole d'Eté de Probabilités de Saint Flour XXXIII. Lecture Notes in Math.* **1896**. Springer, Berlin. To appear.

[17] MASSART, P. and RIO, E. (1998). A uniform Marcinkiewicz–Zygmund strong law of large numbers for empirical processes. In *Festschrift for Miklós Csörgő: Asymptotic Methods in Probability and Statistics* (B. Szyszkowicz, ed.) 199–211. North-Holland, Amsterdam. MR1661481

[18] McDIARMID, C. (1989). On the method of bounded differences. In *Surveys in Combinatorics 1989* (J. Siemons, ed.) 148–188. Cambridge Univ. Press. MR1036755

[19] POLLARD, D. (1982). A central limit theorem for empirical processes. *J. Austral. Math. Soc. Ser. A* **33** 235–248. MR0668445

[20] REYNAUD-BOURET, P. (2003). Adaptive estimation of the intensity of inhomogeneous Poisson processes via concentration inequalities. *Probab. Theory Related Fields* **126** 103–153. MR1981635

[21] TALAGRAND, M. (1996). New concentration inequalities in product spaces. *Invent. Math.* **126** 505–563. MR1419006

[22] TSYBAKOV, A. B. (2004). Optimal aggregation of classifiers in statistical learning. *Ann. Statist.* **32** 135–166. MR2051002

[23] VAPNIK, V. N. (1982). *Estimation of Dependences Based on Empirical Data.* Springer, New York. MR0672244

[24] VAPNIK, V. N. and CHERVONENKIS, A. YA. (1974). *Theory of Pattern Recognition.* Nauka, Moscow. (In Russian.) MR0474638

[25] YANG, Y. and BARRON, A. R. (1999). Information-theoretic determination of minimax rates of convergence. *Ann. Statist.* **27** 1564–1599. MR1742500

[26] YU, B. (1997). Assouad, Fano, and Le Cam. In *Festschrift for Lucien Le Cam: Research Papers in Probability and Statistics* (D. Pollard, E. Torgersen and G. L. Yang, eds.) 423–435. Springer, New York. MR1462963

Equipe Probabilités,
 Statistique et Modélisation
Université de Paris-Sud
Batiment 425
91405 Orsay Cedex
France
E-mail: Pascal.Massart@math.u-psud.fr